\documentclass[a4paper,11pt]{article}
\usepackage[utf8]{inputenc}
\usepackage[english]{babel}
\usepackage[margin=60pt]{geometry}

\usepackage{amsthm,amssymb}
\usepackage{amsmath}
\usepackage{amsfonts}
\usepackage{mathtools} 			
\usepackage[pdfpagelabels]{hyperref}
\usepackage{mathrsfs}
\usepackage{bbm}
\usepackage{stmaryrd}
\usepackage{indentfirst}
\usepackage{thmtools}
\usepackage{thm-restate}
\mathtoolsset{showonlyrefs}

\numberwithin{equation}{section}

\usepackage{bbm}
\usepackage{enumitem}
\usepackage{hyperref}

\usepackage{wrapfig} 
\usepackage{autonum} 

\newtheorem{remark}{Remark}

\usepackage{xcolor}

\newcommand{\matt}[1]{\textcolor{red}{ #1}}
\newcommand{\mv}[1]{\textcolor{blue}{ #1}}

\usepackage{pgf,tikz}
\usetikzlibrary{arrows}
\usepackage{caption}
\usepackage{float} 
\usetikzlibrary{patterns}

\newtheorem{theorem}{Theorem}[section]
\newtheorem{definition}[theorem]{Definition}
\newtheorem{proposition}[theorem]{Proposition}
\newtheorem{corollary}[theorem]{Corollary}
\newtheorem{lemma}[theorem]{Lemma}

\newtheorem{property}[theorem]{Property}

\usepackage{xcolor}

\newcommand{\nahu}[1]{\textcolor{blue}{ #1}}

\newcommand{\NNN}{\mathcal{N}}

\newcommand{\msr}{\mathrm{MSR}}
\newcommand{\vol}{\mathrm{VOL}}
\newcommand{\sr}{\mathrm{Red}}

\newcommand{\pol}{\mathrm{J}}
\newcommand{\Pol}{\{\mathrm{I},\mathrm{II},\mathrm{III}\}}
\newcommand{\MSR}[1]{\msr_{\mathrm{#1}} }
\newcommand{\VOL}[1]{\vol_{\mathrm{#1}} }
\newcommand{\SR}[1]{\sr_{\mathrm{#1}} }
\newcommand{\mr}[1]{\mathrm{#1}}
\newcommand{\NN}{\mathbb N}
\newcommand{\RR}{\mathbb R}

\newcommand{\PPP}{\mathcal P}

\title{Decision-Epoch Matters: 
	Unveiling its Impact on the Stability of Scheduling with Randomly Varying Connectivity}

\author{N. Soprano-Loto$^{a}$, U. Ayesta$^{b,c,d}$, M. Jonckheere$^{a,b}$, and I.M. Verloop$^{b}$ \vspace{10pt}\\
\small
$^a$LAAS-CNRS, Université de Toulouse, CNRS, Toulouse, France\\	\small
$^b$IRIT, Université de Toulouse, CNRS, Toulouse INP, UT3, Toulouse, France\\ \small
$^c$IKERBASQUE, Basque Foundation for Science, 48011 Bilbao, Spain\\ \small
$^d$UPV/EHU, Univ. of the Basque Country, 20018 Donostia, Spain
\normalsize
}

\begin{document}

\maketitle

\begin{abstract}
	A classical queuing theory result states that in a parallel-queue single-server model, the maximum stability region 
	does not depend on the scheduling decision epochs, and in particular 
	is the same for preemptive and non-preemptive systems. 
	We consider here
	the case in which each of the queues may be \emph{connected} to the server or \emph{not}, depending on an exogenous process.
	In our main result, we show that the maximum stability region now \emph{does} strongly depend on how the decision epochs are defined.
	We compare the setting where decisions can be made at any moment in time (the unconstrained setting), to two other settings:  decisions are taken either \emph{(i)} at moments of a departure (non-preemptive scheduling), or \emph{(ii)}  when an exponentially clock rings with rate $\gamma$. 
	We characterise the maximum stability region for the two constrained configurations, allowing us to observe a   reduction compared to the unconstrained configuration. 
	In the non-preemptive setting, the maximum stability region is drastically reduced compared to the unconstrained setting and we conclude that a non-preemptive scheduler cannot take opportunistically advantage (in terms of stability) of the random varying connectivity.
	Instead, for the $\gamma$ decision epochs,  we observe that the maximum stability region is monotone in the rate of the decision moments $\gamma$, and that  one can be arbitrarily close to the maximum stability region in the unconstrained setting if we choose  $\gamma$ large enough. 
	We further show that Serve Longest Connected (SLC) queue is maximum stable in both constrained settings, within the set of policies that select a queue among the connected ones.
	From a methodological viewpoint, we introduce a novel theoretical tool termed a ``test for fluid limits'' (TFL) that might be of independent interest. TFL is a simple test that, if satisfied by the fluid limit, allows us to conclude for stability.
	
\end{abstract}


\section{Introduction}
\label{sec:introd}

In many real-world applications, scheduling decisions are made at specific moments rather than in a  continuous fashion. 
In telecommunication networks, scheduling decisions for allocating resources and managing bandwidth occur at specific moments to adapt to changing traffic patterns. Other examples include server scheduling where certain maintenance tasks or resource allocation decisions are made periodically, or in smart grids where allocation of resources occur at specific moments based on factors like demand patterns and availability of renewable energy. We can also mention real-time systems where once a task is started, it needs to run without interruption to meet its timing requirements. 

It seems intuitively clear  that the performance for the end-user will degrade as the interval between the decision epochs increases since the process becomes ``less controllable.''  
To illustrate this, let us consider a canonical multi-class single-server queue where $\lambda_i$ and $\mu_i$, $i=1, \ldots, K$,  denote the arrival and service rates, respectively, and let us compare the system where decisions can be made at any moment in time (preemptive case) versus the system where decisions can only be made upon a departure of a job (non-preemptive case). In both cases, it is known that the $\mu$-priority policy (which chooses to serve the class in the system with highest $\mu_i$) is optimal. Even though the stability condition for both preemptive and non-preemptive systems is the same  ($\sum_{i=1}^K \lambda_i/\mu_i  <1$), standard analysis of priority policies (\cite[Chapter 5]{Kleinrock76_1}) shows that the mean queue length of the optimal $\mu$-rule in the preemptive case is indeed larger than in the preemptive one. 

In this paper,  we show  that constraining the decision epochs 
can in fact not only have an impact in terms of steady-state performance, but also on the maximum stability region, that is, on the set of parameters for which there exists a policy which induces a steady-state. 
We do so by studying a multi-class scheduling problem with time-varying connectivity, a fundamental problem in networking that has been studied in many papers since the pioneering work by Tassiulas~and~Ephremides \cite{tassiulas1993dynamic}. 
Indeed, time-varying connectivity is inherent in wireless systems and satellite communications~\cite{tassiulas1993dynamic,shakkottai2002scheduling,BM2005,ayesta2013}, where systems face challenges in deciding how to share limited transmitters or channels among users while  weather changes make the transmission rates vary significantly. 
It also arises as a modeling framework in applications where the available capacities fluctuates depending on the traffic conditions. 
More recently, random environments have also been studied in the context of scheduling in data centers and cloud computing (see for instance \cite{peng2015random} or 
\cite{khabbaz2015modelling}), involving a significant amount of communication overhead  such as  data transfers, status updates, and other interactions between the client and the cloud-based services, \cite{dinh2013survey}.
Those communication overhead being variable, the resulting system can be modelled as queuing systems where each queue can be either connected or disconnected to the server according to a random environment.



\begin{wrapfigure}{t!}{0.31\textwidth}
	\centering
	\includegraphics[width=0.30\textwidth]{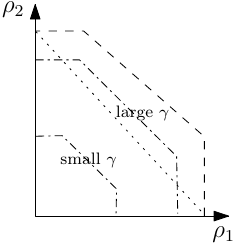} 
	\caption{\small Stability regions with $K=2$ for Settings I (dashed line), II (dotted line) and III (dash-dotted line, for relatively {\em small} and {\em large} $\gamma$). \normalsize }
	\label{fig:your_image_label}
\end{wrapfigure}

We study  a continuous-time system with multiple queues that  can be either connected or disconnected and there is one server that can be dedicated to at most one queue. However, different from previous work, the server is allowed to change to move to another queue  only at certain \emph{specific decision epochs}. We will set out to compare three different settings for these decision epochs. The first setting is when  decisions can be made at any moment in time, which we refer to as the unconstrained setting (as studied in~\cite{BM2005}, referred to as Setting I). The other two are constrained settings where decisions can be taken only at well-defined moments in time:  non-preemptive scheduling (referred to as Setting II), or when an exponentially clock rings with rate $\gamma$ (Setting III). 
In our main result, we characterise the   maximum stability region for the two constrained configurations. The latter are strictly included in the one for  the unconstrained setting, see Figure~\ref{fig:your_image_label}.
To the best of our knowledge  this uncovers a new phenomenon, that is, that the maximum stability region can crucially depends on how decision epochs are defined.
We further show that  policies that select a queue among the connected ones are not necessarily stable due to the random environments, and  prove that the Serve Longest Connected (SLC) queue policy \emph{is} maximum stable for all Settings I, II, and III.

In the non-preemptive case (Setting II),  we obtain that the maximum stability region drastically reduces compared to the preemptive case (Setting I). 
In fact, 
we will show that it coincides with that of a (classical) multi-class single-server queue where the class-$i$ service rate is reduced to the time-averaged departure rate $\mu_i \pi_i(1)$. The latter allows us to infer that in this setting the scheduler cannot take opportunistically advantage (in terms of stability) of the random varying connectivity. 


In Setting III, having decision epochs at an exponential clock with rate $\gamma$  allows us to model applications where there might be a cost associated in observing the state (and taking the corresponding action).
For instance, changing the action too often might be harmful from a performance point of view. 
Therefore, taking decisions with a controllable rate $\gamma$  might be a better option than taking decisions continuously.  
In this respect, our results reveal that the stability condition of Setting~III is strictly smaller than that of Setting~I, but converges to Setting~I  as $\gamma \to \infty$. 
Another interesting example is where a reconfiguration delay is incurred every time a decision is taken, which implies that capacity is wasted. 
The latter means  that the frequency of decision epochs that optimizes the performance is non-trivial.
Our analysis allows us to find the optimal value for the rate $\gamma$ such that the maximum stability region is optimized.


Our main proofs rely on fluid techniques.
The random environment gives rise to complex averaging phenomena, and as a consequence the fluid limits associated to the Markov process under study are in general very complex and cannot be  easily characterized even in low dimensions. On the other hand, the classical quadratic function used  as Lyapunov functions for previous models (either for the stochastic version in the simpler case \cite{tassiulas1993dynamic} or for the fluid limit \cite{shakkottai2002scheduling,andrews2004scheduling}) are not directly applicable to our problem at hand. In our analysis, we introduce a novel theoretical tool termed a ``test for fluid limits'' (TFL), which serves as a simple test that, when satisfied at a fluid scale, allows us to conclude for stability.

Markovian stochastic processes in which actions are taken only at specific decision-epochs, have been  analyzed with the so-called \emph{Embedded Markov Chain} approach. Essentially, this boils down to calculating the transition probabilities and the accumulated reward between two decision epochs, which yields a new Markovian process. This has been the traditional approach to study the performance non-preemptive scheduling in queueing models. 
The novelty of our approach lies in studying
how the stability regions of the underlying Markov processes, which can be seen as the most fundamental performance evaluation criteria, are crucially impacted by the definition of these decision epochs.

\ 

In summary, our main contributions and findings are:
\begin{itemize}
	\item We characterize the maximum stability regions and uncover that they crucially depend on how  decision epochs are defined.
	\item We show that  the SLC policy is maximum stable in all three settings.
	\item 
	Policies that select a queue among the connected ones are not necessarily   maximally stable.
	\item In the non-preemptive setting, we prove that the scheduler cannot take opportunistically advantage (in terms of stability) of the random varying connectivity.
	
	\item In Setting~III, 
	our analysis allows us to find the optimal frequency of decision epochs in order to maximize the  stability region.
	\item We develop a simple test, TFL,  that, if satisfied by the fluid limit, ensures stability. This TFL was used for the three different decision epochs settings.
\end{itemize}


The rest of the paper is organized as follows. In Section \ref{sec:relW}, we recall seminal results from the state of the art for the multi-class scheduling problem with randomly varying connectivity of the queues. The precise model  is described in Section \ref{sec:model}. The maximum stability regions for the three different settings of decision epochs are stated in Section \ref{sec:main}, as well as the different corollaries and insights that follow from them.
In Section~\ref{sec:TestFL} we develop a methodological contribution that proves that a simple TFL criterion allows to conclude for stability. Section~\ref{sec:pmsr} verifies that this TFL criterion holds for our models allowing to prove our characterization of the main stability regions. Finally, numerical experiments and illustrations are reported in Section \ref{sec:num}.

\section{Related work}\label{sec:relW}

As mentioned in the introduction, models with randomly varying connectivity have been studied in the context of wireless and satellite communications, see \cite{tassiulas1993dynamic,shakkottai2002scheduling,Borst05} just to cite a few from a very large body of literature. Other application domains of these models are in data center and  cloud computing, see for instance \cite{peng2015random} or \cite{khabbaz2015modelling}. 

In queueing theory, there is a large body of literature on queues in a random environment. Here, the parameters of the (modulated) queueing model (including arrival rates and service rates) change over time as a function of an exogenous stochastic process. The main focus has been on analyzing the steady-state properties, including stability and queue-length distributions, see for example~\cite{BacceliMakowski}~and~\cite{Dauria}. Our model can be seen as  a particular instance of a modulated queue, where the service rate of queue $i$ fluctuates between 0 and $\mu_i$.

For the parallel-queue single-server model, a seminal paper is \cite{tassiulas1993dynamic} who considered a slotted system with Markovian assumptions and establishes necessary and sufficient conditions for stability. It further showed that SLC is maximum stable and, in the case of symmetric queues, minimizes delay.
The proof technique uses a  quadratic Lyapunov function.
In \cite{BM2005}, a Loyne's construction was employed for an equivalent continuous-time system, wherein the scheduler possesses knowledge of either the workloads or queue-lengths. 
We also refer to \cite{Stolyar04} where it is shown that, under the same stability conditions of \cite{tassiulas1993dynamic}, there exists a so-called Static Service Split rule that is stable.
Similar stability conditions 
were derived in \cite{shakkottai2002scheduling} and \cite{andrews2004scheduling} for more general models (the environment can depend both on classes and the server), and these results were extended to various scheduling types, including variants of the SLC policy and  versions of the Max Weight policy (see the book \cite{srikant2013communication} for a pedagogical review on these models).
See also \cite{celik2011scheduling} for the stability analysis of a model with switching delays.

It is important  to highlight that the fluid limit characterizations derived in \cite{shakkottai2002scheduling} and \cite{andrews2004scheduling}, may not be easily generalized to our context due to the significantly more complex dynamics stemming from the presence of decision epochs during which multiple events can unfold.   
We further note that in all works cited above, the stability results are for models where the  decisions can be made at any moment in time and the state of the environment is 
fully known to the controller. However, to the best of our knowledge,  no stability results have been obtained when  decision epochs are constrained (as in our Setting~II and Setting~III).

One of our main results states that SLC is maximum stable in Settings II and III. Interestingly, SLC is not always maximum stable, as shown for example in \cite{dimakis_walrand_2006} where there are constraints in the set of queues the scheduler can serve. We also refer to \cite{ModianoLQF}, where it is shown that SLC can yield arbitrarily low throughput in a wireless setting with signal interference. 

Finally, let us mention that there are also several works that  deal with an unobservable random environment, that is, when  the decision
maker cannot observe the state of the environment  and it can take its decision only based on
the state of the queues. In~\cite{DURANSIGMETRICS} it is shown that as the number of queues grow large and
the environment changes state relatively fast, an  index policy is asymptotically optimal. Another interesting work is \cite{budhiraja14}, where the authors study a single-server multi-class queueing network in heavy traffic with randomly varying rates. In the main result, it is shown that an ``averaged'' version
of the classical $c\mu$-rule is asymptotically optimal.


\section{Model description}
\label{sec:model}
\subsection{Multi-class scheduling in a random environment}

The system consists of $K$ parallel queues and a single server.
We call $[K]=\{1,\ldots,K\}$ the set of queues.
Each queue $i\in [K]$ has associated an arrival rate $\lambda_i>0$, and a service rate $\mu_i>0$.
For $i\in[K]$, we define the variable $\rho_i=\lambda_i/\mu_i$, usually referred to as the load of the $i$-th queue.
In addition, each queue $i\in [K]$ can be either  connected or disconnected.
Whether or not a queue is connected does not affect the arrivals, but when it is disconnected it cannot receive service.
Formally, the environmental state-space associated to each queue $i\in[K]$ is $\{0,1\}$, the $0$ (resp. the $1$) representing that the $i$-th queue is disconnected (resp. connected).
The environment of each queue has its own law of evolution:
for $ i\in[K] $, the environment of queue~$i$ passes from 0 to 1 (resp. from 1 to 0) at rate $ \lambda_i' $ (resp. at rate $ \mu_i' $).
We denote by $\pi_i$  the invariant distribution of the environment of queue~$i$, which equals
\begin{align}
	\pi_i(0)=\frac{\mu'_i}{\lambda'_i+\mu'_i}
	\quad\mbox{and}\quad
	\pi_i(1)=\frac{\lambda'_i}{\lambda'_i+\mu'_i}.
\end{align}
At each moment in time, the server can be dedicated to at most one queue, the service speed being $\mu_i$ if  queue~$i$ is connected and is being served. 
Decision of which queue to serve can be made only on well-defined \emph{decision epochs}. This decision can not be changed until a next decision epoch. 
There is also an independent of everything else sequence of times, with inter-arrival frequency given by a parameter $ \gamma>0 $,
that can determine the decision epochs.
Decision epochs will be explained in detail in Section~\ref{sec:dec}.
For simplicity of exposition, we will assume that all times between events ---i.e. arrivals, services, changes in environments and the $\gamma$-intensity events just mentioned--- are exponential and mutually independent.

The state space of our  process is
\begin{align}
	\label{lab:ss}
	\mathcal X= \NN_0^K\times \{0,1\}^K\times \{1,\ldots,K,\emptyset\}.
\end{align}
We are adopting the conventions $\mathbb N=\{1,2,\ldots\}$ and $\mathbb N_0=\{0,1,2,\ldots\}$ respectively for the sets of natural numbers and non-negative integers.
The three entries in~\eqref{lab:ss} respectively represent the number of waiting tasks per queue, the states of the environments, and the queue that is receiving service.
We will sometimes refer to the third entry as the state of the server.
If the state of the server is $ \emptyset $, this means that none of the queues receives any service.
A state will be denoted as $ x=(q,e,c) $, and for a given policy and given decision epochs, the process with initial condition $ x $ will be denoted as
\begin{align}\label{eqn:remolacha}
	(X^{x}(t))_{t\ge 0}=((Q^{x}(t),E^{x}(t),C^x(t))_{t\ge 0}.
\end{align} 
We highlight that the distribution of the process depends on the policy applied and the decision epochs considered. However, we chose not to include this information in the notation in order to avoid making it too overloaded, trusting that the context is sufficient to avoid any possible confusion.
At decision epochs $t$, the policy is assumed to make decisions based only on the pair $(Q^x(t),E^x(t))$, namely from the present time observation of the queue sizes and the states of the environment, in particular implying that the process is Markov.

\begin{remark}[Markovian assumptions]
	We opted for Markovian assumptions on the dynamics because of an escalation in technical difficulty without a significant gain in insights. In particular, dealing with more general distributions 
	would give rise to non-explicit and intricate expressions of the key parameters, while their Markovian counterparts are simple.   
\end{remark}


\subsection{Decision epochs}
\label{sec:dec}
We describe now the three different settings for decision epochs: 

\begin{enumerate}[label=\textbf{Setting \Roman*}.,align=left]
	\item Decision epochs are at any moment in time. 
	Stability under this setting has been studied originally in discrete time in \cite{tassiulas1993dynamic}, while continuous-time versions and generalisations have been considered in \cite{bambos2004queueing,shakkottai2002scheduling} (see Section~\ref{sec:relW} on related).
	\item
	This is the non-preemptive situation in which a decision cannot be taken in the middle of a task. In other words, decision epochs are any moment after a departure and until a new task enters into service.
	
	\item Decision epochs are when an exogenous exponential clock of intensity $ \gamma >0$ rings.
	
\end{enumerate}

In the context of queues with random connectivity, to the best of our knowledge, Setting I has been the only case studied so far. We can mention \cite{tassiulas1993dynamic,shakkottai2002scheduling,BM2005,ayesta2013} and other references in Sections~\ref{sec:introd}~and~\ref{sec:relW}. Setting~II is motivated by applications  in which a decision cannot be taken in the middle
of a task. This occurs for instance in real-time systems where once a task is started, it needs to run without interruption to meet its timing requirements. Finally, Setting~III aims at modeling situations in which there might be a cost associated to observing the state or to taking a decision. Here, the administrator might need to define the scheduling epochs at some intervals (defined by clock intensity $\gamma$) in order to strike a good balance between performance and cost (both being increasing as a function of $\gamma$).


\subsection{Scheduling policies and maximum stability regions}

We say that a policy stabilizes the system if its associated Markov process is positive recurrent or, equivalently, if it has a unique invariant distribution.	
For $ \pol\in\Pol  $, let $ \msr_\pol $ be the maximum stability region associated to Setting $ \pol $,
which is  defined as follows:
\begin{align}
	\msr_\pol=\{ (\lambda,\mu,\lambda',\mu'):\mbox{there exists a policy that stabilizes the system under Setting \  $\pol$}  \}.
\end{align}
We are using the compact notation $\lambda=(\lambda_1,\ldots,\lambda_K)$, $\mu=(\mu_1,\ldots,\mu_K)$, and so on.
A policy is called \emph{maximum stable} if it can stabilize the system for any set of parameters $(\lambda,\mu,\lambda',\mu')\in \msr_\pol$. 

Characterizing the maximum stability region for Setting III, $\MSR{III}$, seems out of reach. We refer to Section~\ref{sec:SLC_is_not_MS} for a full discussion on this. In order to obtain a closed-form expression for the maximum stability region, we will restrict the search to the following class of policies, denoted by~$\mathcal{P}$. For insights as to why we restrict to this set, we refer to the text right above Section~\ref{sec:maxstable}.
\begin{definition}[class $\PPP$ of policies]\label{def:admissible}
	We say that a policy is inside the class $\mathcal P$ if at decision epochs it can only choose to serve a queue that is connected. In addition, it will choose a connected queue with waiting tasks, whenever possible.
	The policies in this family will be called $\PPP$-\textit{policies}.
\end{definition}

Note that when no queue is connected at a decision epoch, no queue is served and hence the state of the server is $\emptyset$, until the next decision epoch.
An important policy inside the class~$\PPP$ is the Serve the Longest Connected (SLC) policy, as first defined in \cite{tassiulas1993dynamic}. SLC is defined as the policy that \emph{at each decision epoch} chooses to dedicate the server to the queue with the highest number of waiting tasks among the ones that are connected, with an arbitrary tie-breaking rule.

The maximum stability region under Setting $\pol\in\Pol$ when restricted to policies in the family $\PPP$ is defined as:
\begin{align}
	\msr_\pol^{\mathcal P}=\{ (\lambda,\mu,\lambda',\mu'):\mbox{there exists a   $\mathcal P$-policy that stabilizes the system under Setting $\pol$}  \}.
\end{align}
We will write $\MSR{III}(\gamma)$ to show the dependence on $\gamma$ in the case  of Setting~III.

\section{Main stability results}
\label{sec:main}
In this section, we state our main result, which gives an explicit characterization of the maximum stability regions in all three settings.

\begin{theorem}
	\label{thm:MSR}
	For $L\subseteq[K]$, let 
	\begin{align}
		\pi_L^{0}=\prod_{i\in L}\pi_i(0)
	\end{align}
	represent the probability that all the environments of queues in $L$ are disconnected.
	We then have that
	\begin{enumerate}
		\item $\MSR{I}$ and $\MSR{I}^\mathcal{P}$ are characterized  by
		\begin{align}
			\label{eqMSRI}
			\sum_{i\in L}\rho_i<1-\pi_L^{0}\quad \forall  L\subseteq [K], L\neq\emptyset.
		\end{align}
		\item $\MSR{II}$ and $\MSR{II}^\mathcal{P}$ are  characterized by
		\begin{align}\label{eq:condition_II}
			\sum_{i=1}^K\frac{\rho_i}{\pi_i(1)}<1.
		\end{align}
		\item Finally,  $\MSR{III}^{\mathcal P}(\gamma)$ is characterized by
		\begin{align}
			\label{eqMSRIII}
			\sum_{i\in L}\frac{\rho_i}{\theta_i(\gamma)}<1-\pi_L^{0}\quad \forall L \subseteq [K],L\neq\emptyset,
		\end{align}
		where
		\begin{align}
			\theta_i(\gamma)=\frac{\gamma+\lambda_i'}{\gamma+\lambda_i'+\mu_i'}.
		\end{align}
	\end{enumerate}
	In all the three settings, the corresponding maximal stability regions are attained by the SLC policy.
\end{theorem}

For Setting~I, both the maximum stability condition and the stability of SLC   were first obtained in~\cite{tassiulas1993dynamic} (see also \cite{BM2005}), while the results for Settings~II and III are new.
The proof of Theorem~\ref{thm:MSR} can be found in Section~\ref{sec:pmsr}. We note that the proof (for all three settings) is based on a novel methodological tool termed TFL, which serves as a simple test to conclude for stability, see Section~\ref{sec:TestFL}. 

Intuitively, the   maximum stability regions under the three settings can be explained as follows.
Consider a subset of classes~$L$.
The right hand side (RHS) represents the maximum fraction of time the server can be  usefully dedicated to this set~$L$ in a saturated regime and when queues in $L$ are given priority over the other queues.
The left hand side (LHS) of the conditions represents the \emph{effective load} of this subset~$L$, which is defined as $\sum_{i\in L} \lambda_i/(\xi_i^{\,\pol}\mu_i)$, $\pol\in \Pol$.
Here, $\xi_i^{\,\pol}$ represents the proportion of time the $i$-th queue is connected within the time the server is usefully dedicated to it.
In this way, $\xi_i^{\,\pol}\mu_i$ represents the  effective departure rate.
The proportion $\xi_i^{\,\pol}$ appears because, due to the restrictions on decision making,
there is a proportion of time in which the server is stuck in a disconnected queue,
even though other queues in $L$ are connected.
This phenomena does not occur in Setting I, reason why $\xi_i^{\,\rm I}=1$.
The three conditions now follow: 
\begin{itemize}
	\item
	Setting I: The RHS is $1-\pi^0_L$, since it is only useful to dedicate the server to the set $L$ when at least one queue in $L$ is connected.
	As already mentioned, 
	since decision epochs are any moment in time,
	the server only dedicates service to a queue when it is connected, so that  the effective departure rate of a class is $\mu_i$,
	which explains the LHS.
	\item
	Setting II: Since the service is non-preemptive, we are sure that at a decision epoch the queue under service is connected. We can therefore stay serving a queue in the set~$L$, since for sure there was one queue in~$L$ connected at the decision epoch. 
	The fraction of time the server can be dedicated to a subset~$L$  is therefore simply 1, that is, the RHS.
	We briefly explain how the quantity $\xi_i^{\,\rm{II}}=\pi_i(1)$ can be obtained.
	For a fixed queue $i$, the process encoding the state of its environment (connected or disconnected) only when the server is dedicated to this queue, 
	turns out to be a Markov process
	whose invariant distribution is precisely $\pi_i$.
	Then the ergodic theorem implies that the desired proportion $\xi_i^{\,\rm{II}}$ is $\pi_i(1)$.
	A more detailed explanation is provided in Section~\ref{sec:TFLtwo}.
	
	The reason why the maximum stability condition does not depend on subsets of queues is that,
	if the stability condition is satisfied for the total queue subset $[K]$,
	then it is automatically fulfilled for any other subset $L$, since  the RHS does not depend on the subset of queues $L$ considered.
	\item
	Setting III:  At a decision moment, the server can be dedicated to a queue in~$L$ only if one of the queues is connected (recall that we are restricting to $\PPP$-policies). 
	The fraction of time in which this occurs is $1-\pi^0_L$ because
	the law governing the decision epochs (the exponential clocks of rate $\gamma$) is independent of the law governing the environments.
	This explains the RHS.
	The reason why $\xi_i^{\,\rm{III}}=\theta_i(\gamma)$ is similar to that in Setting II,
	in this case the Markov process at issue having an invariant distribution that depends on $\gamma$.
	A detailed explanation can be found in Section \ref{sec:TFLthree}.
\end{itemize}

\begin{center}
	\begin{figure}[t]
		\includegraphics[width=1\textwidth]{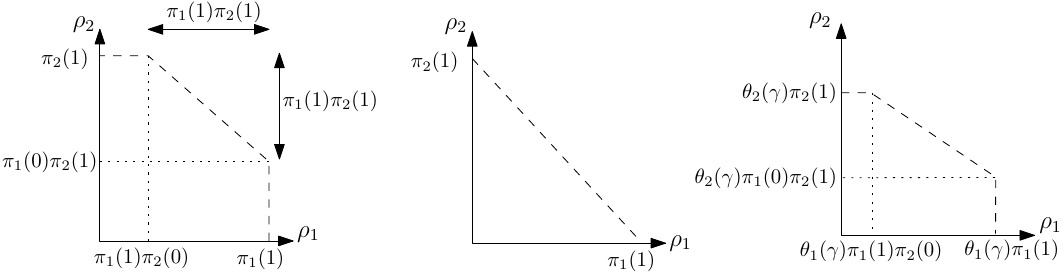}
		\caption{For $K=2$, we depict $\MSR{I}$ {\em (left)}, $\MSR{II}$ {\em (center)} and $\MSR{III}$ {\em (right)}.
		}
		\label{fig:stabilityregion}
	\end{figure}
\end{center}

Our result shows that  the maximum stability region strongly depends on the chosen setting of the decision epochs.  
In particular, the $\MSR{\mr{II}}$  and $\MSR{\mr{III}}^{\mathcal{P}}(\gamma)$ are strict subsets of $\MSR{\mr{I}}$ as expected, since the decision epochs are more restricted. 
We refer to Figure~\ref{fig:your_image_label} for a plot of the different maximum stability regions for the case $K=2$.




There is a small caveat concerning Setting III. In fact, the maximum stability region is obtained when restricting to  policies in $\PPP$ only. In Section~\ref{sec:SLC_is_not_MS}, we will prove that policies outside this class can outperform the $\PPP$-policies in terms of stability. That is, from the stability point of view, it might be better to serve a queue that is currently disconnected, even though there is a connected queue with work waiting. This counter-intuitive property results from the fact that after a decision is made, disconnected queues can become connected (before the next decision epoch), and  hence the occasional strategic choice of serving a currently disconnected queue can therefore optimize for higher departure rates in the future. Hence, the actions taken under a maximum stable policy should be a function 
of whether or not queues are connected which complicates significantly the analysis. 
Due to this complexity, we leave the search for such a maximum stable policy, as well as finding the maximum stability region $\MSR{III}(\gamma)$, for future research and focus here on $\MSR{III}^\PPP(\gamma)$ instead.


We note that restricting to the set of policies  $\PPP$ is relevant. Firstly, this allows us to obtain interesting results such as the characterization in closed-form of the maximum stability region and in showing that SLC is maximum stable. Secondly, our results show that $\mathcal{P}$ is somehow a ``complete'' class of policies, since by taking $\gamma$ large enough its stability region can be arbitrarily close to the maximum stability region in the unconstrained setting (Setting~I).


\subsection{Are  $\PPP$-policies always maximum stable?}
\label{sec:maxstable}

In the  classical multi-class scheduling problem, that is, when all the queues are always connected, and decisions can be taken at any moment in time, it is well known that any \emph{work-conserving policy} is stable in the maximum stability region. In fact, this is valid for both preemptive and non-preemptive systems.
A natural question is whether this remains valid when adding either \textit{(a)} varying connectivity or \textit{(b)} decision restrictions defined by~$\gamma$. We note that in the presence of \textit{(a)} and/or \textit{(a)}, the set of   policies $\PPP$ is the equivalent of work-conserving policies in the classical queue. 
Thus, in this section we aim at answering the following question, {\bf \em  Q1: Is it true that any $\PPP$-policy is maximum stable  in the presence of \textit{(a)} or \textit{(b)}? }
We note that Setting~I and Setting~II would correspond to adding \textit{(a)} to a preemptive and non-preemptive system, respectively, while Setting~III would correspond to having both \textit{(a)} and \textit{(b)}.
When instead we add only~\textit{(b)}, we retrieve a classical multi-class scheduling problem where decisions can be taken only at $\gamma$-epochs, which we refer to as Setting~0 (with slight abuse of notation).
In the proposition below, we show that the answer to \emph{\bf Q1} is negative for all cases. The underlying reason being that, unlike in the classical case, by adding  \textit{(a)} and/or \textit{(b)} to the system, for stability it is important to keep queues balanced. That is, it is risky to have queues with a small number of tasks, because  in the case of \textit{(a)} this leads to states where all connected queues are empty and we have to idle unnecessarily, or  in the case of \textit{(b)} deciding at a $\gamma$-moment to serve a connected queue that has very little tasks brings the risk that the queue empties before the next decision moment and as such the server is again unnecessarily idle.


\begin{proposition}
	\label{prop:prio}
	Given a decision epoch Setting~$\pol$, with $\pol=0, \rm{I}, \rm{II}$ or $\rm{III}$, and given a priority ordering of the queues.
	We consider the priority policy that at each decision moment serves the queue with highest priority among all connected queues with waiting tasks.
	Then, there exist parameters inside the maximum stability region $\MSR{J}^\PPP$ for which this priority policy is not stable. 
\end{proposition}

This provides the negative answer to \emph{\bf Q1} because the priority policy defined in Proposition~\ref{prop:prio}  is inside the class~$\PPP$, but it is not maximum stable. 
Below we give the proof for $J=I$,  the other settings follow in a similar fashion.

\begin{proof}
	Assume we are in Setting~I and w.l.o.g.\ assume we have the  priority ordering $1\succ 2\succ \ldots \succ K$.  
	We consider $\pi_i$, $i=1,\ldots,K$, fixed and take $\rho_1 \approx 0$.
	Since  $\rho_1<1-\pi_1(0)$,  the dynamic of the sole queue $1$ is stable,
	so that we can define the asymptotic proportion of time~$\varepsilon_1>0$ that this queue is non-empty and connected.
	Since queue~1 has priority, queue~2 can only be served when queue~1 is not served, which happens with a fraction of time~$1-\varepsilon_1.$
	Queue~2 can use this time only if it is connected, hence, it is stable if and only if 
	\begin{equation}
		\label{eq:vareps}
		\rho_2 < (1-\pi_2(0))(1-\varepsilon_1).
	\end{equation}
	Observe that $\varepsilon_1$ has been constructed independently of $\rho_2$, so at this point we can chose $\rho_2$ in such a way that  condition~\eqref{eqMSRI} is fulfilled for $L=\{2\}$ and $L=\{1,2\}$, while~\eqref{eq:vareps} is not. Now setting also $\rho_3=\ldots =\rho_K\approx 0,$ we can guarantee that all conditions in~\eqref{eqMSRI} are fulfilled.
	Hence, even though the parameters are inside the maximum stability region~$\MSR{I}$, the priority policy is not stable.
\end{proof}

\subsection{Non-preemptive versus preemptive scheduling}
\label{sec:p-vs-np}
In the non-preemptive case (Setting II), Theorem~\ref{thm:MSR} shows that the maximum stability region drastically reduces compared to the preemptive case (Setting I). As explained previously, in Setting~II, the stability region coincides with that of a (classical) multi-class single-server queue in which the class $i$ service rate is  $\mu_i \pi_i(1)$. 

The reduction in stability due to the non-preemption assumption can be calculated as follows:
From Theorem~\ref{thm:MSR} it follows that the maximum stability region  is a convex-hull whose volume can be computed. 
In order to assess the degradation of the stability region we compare    $\VOL{I}$ and $\VOL{II}$, where $\VOL{I}$ ($\VOL{II}$) denotes the volume of the space enclosed by points $(\rho_1,\cdots,\rho_K)$ that satisfy $\MSR{I}$ ($\MSR{II}$). 
In the  case $K=2$, see Figure~\ref{fig:stabilityregion}, we can calculate these volumes: $\VOL{I}= \pi_1(1) \pi_2(1) - \frac{(\pi_1(1)\pi_2(1))^2}{2}$ and  $\VOL{II}= \frac{ \pi_1(1) \pi_2(1)}{2}$, which gives as reduction in stability 
\begin{equation}
	\label{eq:redII}
	\SR{II} := \frac{\VOL{I}- \VOL{II}}{\VOL{I}} 
	=  \frac{1-\pi_1(1)\pi_2(1)}{2- \pi_1(1)\pi_2(1)}, 
\end{equation}
for $\pi_1(1), \pi_2(1)>0,$ due to the non-preemptive assumption. 
This shows that the reduction does not depend on the arrival rates or departure rates. It does however depend on the stationary distribution of the connectivity of the queues. We see that the reduction is close to zero when the queues are almost always connected, while it is close to 1/2 when the queues are almost always disconnected. 
To illustrate the latter, consider Figure~\ref{fig:stabilityregion} and observe that $\MSR{I}$ converges to a square (whose volume converges to zero), and the border of $\MSR{II}$ is the diagonal of this square. Hence, in the limit we obtain as reduction 1/2.


For arbitrary $K$, $\MSR{I}$ and $\MSR{II}$ are polytopes defined by the convex hull of the vertices obtained from  linear inequalities. We  note that there is no analytical expression to calculate the volume of a convex polytope in general. However, $\MSR{II}$ has a special structure, a simplex, for which the general formula for the volume is $\VOL{II} 
=\frac{1}{K!}\prod_{i=1}^K\mu_i \pi_i(1).$ 
Unfortunately, since we do not have any expression for $\VOL{I}$, $\SR{II}$ can only be estimated by numerical means. 

\subsection{Different time scale regimes}
\label{sec:dts}
In this section, we study the impact that the speed of the environment has over the stability regions.
For that purpose, we use a speed factor $\alpha>0$
which acts over the system by multiplying the rates $\lambda'_i$ and $\mu'_i$ for each $i$.
The first observation is that $\msr_{\rm{I}}$ and $\msr_{\rm{II}}$ remains unaffected under the scaling given by $\alpha$.
In Setting III, where this is not the case anymore, we will study
the impact that this scaling has in combination with the parameter $\gamma$.



\begin{itemize}
	\item If $\frac{\alpha}{\gamma} \to 0$, then from Theorem~\ref{thm:MSR} it follows that 
	\begin{equation}
		\label{eq:limit}
		\lim_{\frac{\alpha}{\gamma} \to 0} \MSR{III}^\PPP (\gamma) \to \MSR{I}.
	\end{equation}
	This can be explained because, in this limiting regime,
	decision epochs are so often compared to changes in the environment that we approach Setting I in which decisions can be made at every moment.
	
	
	\item If $\frac{\alpha}{\gamma} \to \infty$, then from Theorem~\ref{thm:MSR} it follows that  
	\begin{align}
		\lim_{\frac{\alpha}{\gamma}\to\infty} \MSR{III}^\PPP(\gamma)  \to     \sum_{i\in L}\frac{\rho_i}{\pi_i(1)}<1-\pi_L^{0}\quad \forall L \subseteq [K],L\neq\emptyset.
	\end{align}
	It can be easily checked that the latter is strictly included in  $\MSR{II}$.
	To explain this inclusion, observe that, 
	in the limiting regime,
	the states of connectivity of the queues changes infinitely many times between two decision epochs, and as such one observes as effective departure rate the steady state $\mu_i \pi_i(1).$
	However, contrary to Setting~II, the fraction of time that can be dedicated to a set of queues~$L$ is given by $1-\pi_L^{0}$, since only queues that are connected at a $\gamma$ moment are allowed to receive service, which explains the strict reduction in stability in the RHS.  
\end{itemize}



\subsection{Sufficiently large   frequency to guarantee stability}
\label{sec:III-vs-II}

It follows from Theorem~\ref{thm:MSR} that $\MSR{III}^\PPP(\gamma)$ is monotone increasing on $\gamma$, and in particular, $\MSR{III}^\PPP(\gamma) \to \MSR{I}$ as $\gamma \to \infty$, see Section~\ref{sec:dts}. In some applications, there might be a cost (for instance energy consumption) associated to observing the state. In such a situation, the system administrator might want to keep  the   frequency $\gamma$, as low as possible, while keeping the system stable. 
Theorem~\ref{thm:MSR} ensures the existence of a minimal value for  $\gamma$, denoted by $\gamma_0$,  such that the system can be made stable for any $\gamma>\gamma_0.$ 
\begin{corollary}
\label{cor:gamma0}
For a set of parameters $(\lambda,\mu,\lambda',\mu')\in \msr_I$, there exists a finite $\gamma_0$ such that for $\gamma> \gamma_0$, the SLC policy is stable in Setting~III.
\end{corollary}

Equivalently to  Subsection~\ref{sec:p-vs-np}, we define  $\VOL{III}(\gamma)$ as the volume of the space enclosed by points $(\rho_1,\cdots,\rho_K)$ that satisfy $\MSR{III}^\PPP(\gamma)$, and $\SR{III}(\gamma) := \frac{\VOL{I}- \VOL{III}(\gamma)}{\VOL{I}}$ as the reduction of stability for Setting~III compared to Setting~I. 
Even though we cannot calculate the volumes $\MSR{I}$ and $\MSR{III}^\PPP(\gamma)$ for $K>2$, it turns out we can instead directly calculate $\SR{III}(\gamma)$. 
To see this, we note that if $(x_1,x_2,\ldots,x_K)$ is a vertex of $\MSR{I}$, then $(\theta_1(\gamma) x_1, \theta_2(\gamma) x_2,\ldots, \theta_K(\gamma) x_K)$ is a vertex of $\MSR{III}^\PPP(\gamma)$ as well. 
We refer to Figure~\ref{fig:stabilityregion} for an illustration when $K=2$.
This then allows us to relate their volumes, see the following corollary. For the proof, we refer to the Appendix~\ref{sec:addproofs}.

\begin{corollary}
\label{eq:corredIII}
For arbitrary $K$, we have 
\begin{equation}
	\label{eq:redIII}
	\SR{III}(\gamma):= \frac{\VOL{I}- \VOL{III}(\gamma)}{\VOL{I}} =  1 - \theta_1(\gamma) \theta_2(\gamma) \cdots  \theta_K(\gamma). 
\end{equation}
\end{corollary}


Since $\theta_i(\gamma)$ is increasing in $\gamma$, it follows that  $\SR{III}(\gamma)$ decreases as $\gamma$ increases. We can then consider a similar example to the one that led to  Corollary~\ref{cor:gamma0}. Let us assume that   there is cost associated to observing the state. A system administrator might then want to determine the rate $\gamma_1$ sufficiently large so that the stability reduction is smaller than a certain threshold. 

\begin{corollary} 
\label{cor:gam1}
For a set of parameters $(\lambda,\mu,\lambda',\mu')\in \msr_I$ and a target stability reduction $R$, there exists a finite $\gamma_1$ such that $\SR{III}(\gamma) \leq R$, for all $\gamma \geq \gamma_1$.
\end{corollary}

\subsection{
Optimizing the maximum stability region in the presence of communication overhead}
\label{sec:commoverhead}

Assume that each decision-making incurs a communication overhead, leading to the server being momentarily diverted from task execution and engaged in a communication process.
To make things simpler, let us further assume that each decision-making then incurs an inactivity of fixed duration $c$ before the server resumes. 
For this model, the maximum stability region can be derived from $\msr_{\rm{III}}^{\mathcal{P}}(\gamma)$ and is given in the following Corollary (see the Appendix~\ref{sec:addproofs} for the proof).
\begin{corollary}\label{optigamma}
Consider Setting $\rm{III}$. Assume that at each decision epoch the server incurs an inactivity period of fixed duration~$c$. 
The precise scheduling decision of which queue to serve is decided after this inactivity period.
The maximum stability region is then given by
\begin{align}
	\label{eq:optrateoverheads}
	\sum_{i\in L}\frac{\rho_i}{\theta_i(\gamma)}< { 1-\pi_L^{0} \over 1+\gamma c}\quad \forall L \subseteq [K],L\neq\emptyset.
\end{align}
It then follows that for a given set of parameters, there exists a finite $\gamma_c^*$ such that the maximum stability region is maximized in $\gamma^*_c$.
\end{corollary}

The value of $\gamma_c^*$ corresponds to the optimal rate of decision epochs that strikes the right balance between a high frequency of taking actions and the overheads induced by them.

\subsection{Could SLC  be improved upon?}
\label{sec:SLC_is_not_MS}

The SLC policy was shown to be maximum stable when decisions can be made at any moment in time (Setting~I), as well as when decisions have to be non-preemptive (Setting~II), see Theorem~\ref{thm:MSR}.
When decisions happen at $\gamma$-moments (Setting~III), we proved that SLC is maximum stable when restricting to the set of $\PPP$-policies. The question remains whether SLC is also maximum stable in general, that is, stable in the set $\MSR{III}(\gamma)$. In  this section, we show that the answer to this question is negative and discuss other policies that could potentially be maximum stable. 

To show that SCL is not maximum stable  in Setting III, we  prove  that 
the maximum stability region $\MSR{III}(\gamma)$ is strictly larger than $\MSR{III}^\PPP(\gamma)$. This is stated in Proposition~\ref{prop:strict} and for its proof we refer to Appendix~\ref{sec:addproofs}.
Intuitively, the first inclusion in~\eqref{eq:aa} can be seen as follows: 
we will define a policy outside the family $\PPP$
which improves SLC in terms of stability. That is, there exist sets of parameters outside $\MSR{III}^\PPP(\gamma)$ for which this new policy is stable.
The policy is defined as follows: at a decision epoch, its action is identical to SLC if there is at least one  queue that is connected; otherwise, when all the queues are disconnected, it dedicates the service to the queue with the highest number of tasks instead of going to the $\emptyset$ state as SLC does. The improvement in terms of stability comes from the fact that the intervals delimited by consecutive $\gamma$-marks at which SLC is in state $\emptyset$ represent a positive proportion of time which is seized by the new policy.

\begin{proposition}
\label{prop:strict}
It holds that
\begin{equation}
	\label{eq:aa}
	\MSR{III}^\PPP(\gamma) \subset \MSR{III}(\gamma)  \subset \MSR{I}.
\end{equation}
Here $\subset$ denotes a  subset in the strict sense.
\end{proposition}

The second inclusion in the above proposition states that   $ \MSR{III}(\gamma)\subset \MSR{I}$. That is, the timing restriction in  Setting III creates a strict difference in terms of stability.
Nevertheless, as stated in Corollary~\ref{cor:gamma0}, in Setting~III there exist policies (e.g. SLC) that are  stable for parameters in $\MSR{I}$ as the decision rate~$\gamma$ is set sufficiently large.
A full characterization of $\MSR{III}(\gamma)$ is left as future research. 

Another question that is left unanswered is which policies \emph{do} stabilize the system for parameters in $\MSR{III}(\gamma).$
For now, we do not have an answer to this. We do  believe however that \emph{a maximum stable policy will need to occasionally serve a disconnected queue} even though there are connected queues with tasks waiting to be served\footnote{Note that the earlier mentioned ``new version'' of SLC was only different to SLC  when all queues where disconnected. This was sufficient to prove the strict inclusion of $\MSR{III}^\PPP$ in $\MSR{III}$, but it will not be enough  for the policy to be maximum stable.}:
Although choosing a disconnected queue is disadvantageous in the short term, it might be beneficial in the long run. The latter  follows because between two decision epochs the queue might become connected and hence be able to serve tasks. If the disconnected queue had a large backlog while the connected queue had very few tasks waiting, serving the disconnected queue might be the preferred action  in order to avoid unbalanced queues (see  Section~\ref{sec:maxstable} for an explanation as to why one would want to avoid unbalanced queues). 

We conclude this section by mentioning two type of policies that might be able to provide a maximum stable system in Setting~III.
The first one is a static policy (that is, whose actions do not depend on the queue lengths) defined by a  so-called Static Service Split (SSS) rule  as introduced in~\cite{Stolyar04} for a more general environment model (without restrictions on the decision epochs). Under the SSS policy, each possible set of connected queues is associated a probability vector (the static service split) that determines with which probability each of the queues is chosen to be served at a decision epoch when only this set of queues is connected. 
For Setting~I, stability of SSS rules was proved in~\cite[Proposition 1]{Stolyar04}. More precisely, it is proved that for each set of parameters in the maximum stability region~$\MSR{I}$ there exists a static service split such that the SSS rule is stable in Setting~$I$.
We believe that similar results hold true for $\MSR{III}(\gamma)$ in Setting~III, but leave this for future research.
Another candidate for a stable policy we would like to mention is based on the restless bandit framework as introduced by Whittle in \cite{Whi88}. In this framework, a queue~$i$ is seen as an arm whose two-dimensional state consists of the number of waiting tasks, $Q_i(t)$,  and whether or not the queue is connected, $E_i(t)$. Each queue is associated a Whittle index as a function of $Q_i(t)$ and $E_i(t)$. The so-called Whittle index policy consists in serving at each decision epoch the queue that has currently the highest Whittle index value. Results in the literature show that such a policy can be very efficient, see~\cite{weber1990index,Verloop16}, and we plan to further investigate its stability properties in the context of our model.

\section{Test for fluid limits (TFL) and stability}
\label{sec:TestFL}

This section is devoted to introducing a new methodological strategy to prove stability, that we call \textit{test for fluid limits} (TFL).
In Section~\ref{sec:pmsr} we use this method to prove all the settings considered in the article.
Moreover, we believe that this strategy might be applied to cover other contexts.
The idea behind this approach is to obtain the stability of the system as a consequence of the verification of a certain test for the fluid limit, 
the advantage of course being that this verification is easy to carry out.
In this way, we succeed in proving stability without the necessity of explicitly describing the fluid limits, which turns out to be cumbersome in our examples.
Before stating Proposition \ref{lemma:test}, the result enclosing this idea, certain prerequisites must be defined.

We introduce the notion of fluid limit following \cite{Bra2008}.
All our Poisson point processes used to construct our stochastic process are defined in the same probability space $(\Omega, \mathcal F, \mathbb P)$.
We fix an almost sure event $\Omega'\in \mathcal F$, 
which can be defined as a set for which the law of large numbers holds for all these Poisson processes.
Let $\|\cdot\|$ be the supremum norm in $\RR^K$.
For $x=(q,e,c)\in \mathcal X$ such that $\|q\|\neq 0$, we define the random function $\bar Q^x:[0,\infty)\to\mathbb R^K$ as
\begin{align}\label{scaledQ}
\bar Q^x(t)=\frac{1}{\|q\|}Q^x(\|q\|t).
\end{align}
This function is the rescaled queue length process.
A sequence $(x^n)_n=((q^n,e^n,c^n))_n\subseteq \mathcal X$ is said to be divergent if $\lim_{n\to\infty}\|q^n\|=\infty$. 

\begin{definition}
A \textit{fluid limit} is a (deterministic) function $G:[0,\infty)\to \mathbb R^K$ for which there exist $\omega\in\Omega'$ and a diverging sequence $(x^n)_n$ such that
\begin{align}\label{eq:album}
\lim_{n\to\infty}\bar Q^{x^n}(t)[\omega]=G(t),
\end{align}
where the convergence is component-wise and uniformly over compact time subsets.    
\end{definition}



To a function $ G=(G_i)_i:[0,\infty)\to \RR^K $,
we associate the max and the argmax functions, which are defined as
\begin{align}
M(t)=\max_{i\in[K]}G_i(t)\quad \text{and}\quad  L(t)=\{i\in[K]:G_i(t)=M(t)\}, \quad t\ge 0.
\end{align}
The definition of our test for fluid limits (TFL) follows.

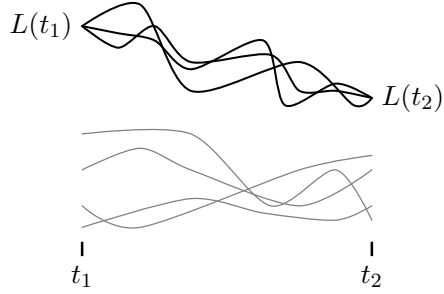
\begin{figure}[htbp]
\centering
\resizebox{!}{4cm}{
\begin{tikzpicture}
	\draw[thick,line cap=round] plot [smooth, tension=0.5] coordinates {(-2,1) (-1.5,0.7) (-1,1) (-0.5,0.5) (0.6,0.6) (1,0.1) (1.5,0.1) (2,0)};
	\draw[thick,line cap=round] plot [smooth, tension=0.5] coordinates {(-2,1) (-1.5,0.9) (-1,0.8) (-0.5,0.4)  (0.5,0.8) (0.8,-0.1) (1.5,0.2) (2,0)};
	\draw[thick,line cap=round] plot [smooth, tension=0.5] coordinates {(-2,1) (-1.2,1.3) (-0.5,0.1)   (1,0.5)  (1.75,-0.1) (2,0)};
	\draw[gray,line width=0.5pt,line cap=round] plot [smooth, tension=0.5] coordinates {(-2,-0.5)  (-0.5,-0.5) (0.6,-1.5) (1.5,-1) (2,-1.7)};
	\draw[gray,line width=0.5pt,line cap=round] plot [smooth, tension=0.5] coordinates {(-2,-1.8) (-0.5,-1.4)  (0.5,-1.6)  (1.5,-1.7) (2,-1.5)};
	\draw[gray,line width=0.5pt,line cap=round] plot [smooth, tension=0.5] coordinates {(-2,-1) (-1.2,-0.7) (-0.5,-1)   (1,-1.5)  (2,-1)};
	\draw[gray,line width=0.5pt,line cap=round] plot [smooth, tension=0.5] coordinates {(-2,-1.5) (-1.2,-1.8)   (1,-1)  (2,-0.8)};
	\node[left] at (-2,1) {$L(t_1)$};
	\node[right] at (2,0) {$L(t_2)$};
	\draw[line width=1pt] (-2,-2) -- (-2,-2.2);
	\node[below] at (-2,-2.2) {$t_1$};
	\draw[line width=1pt] (2,-2) -- (2,-2.2);
	\node[below] at (2,-2.2) {$t_2$};
\end{tikzpicture}
}
\caption{This picture represents condition \eqref{eqn:ladrillo}.
The three leading lines are the graphs of the functions $G_i$ for which $i$ attains the maximum in the extremes of the interval under consideration, $[t_1,t_2]$.
The other lines are the graphs of the remaining $G_i$'s.
The two groups of functions are not allowed to intersect among the hole interval $[t_1,t_2]$.}
\label{fig:property1}
\end{figure}

\begin{definition}[Test for fluid limits (TFL)]\label{def:property_1}
We say that the TFL is passed if there
exists $ \delta>0 $ such that
for every fluid limit $G$
and every pair of times
$ 0\le t_1<t_2 $ satisfying
\begin{align}
\label{eqn:ladrillo}
L(t_1)=L(t_2)
\text{ and }
\displaystyle{\min_{i\in L(t_1)}G_i(t)}>\max_{i\notin L(t_1)}G_i(t)\text{ for every }t\in [t_1,t_2],
\end{align}
we have
\begin{align}\label{eq:verdura}
\frac{M(t_2)-M(t_1)}{t_2-t_1}\le -\delta.
\end{align}
In Figure~\ref{fig:property1},
we illustrate condition \eqref{eqn:ladrillo}. 
\end{definition}

The following result is our methodological tool, that guarantees stability once the test is passed.
Its proof is given in Section~\ref{sec:proof_test} and relies on an inductive procedure that, to the best of our knowledge, is novel.

\begin{proposition}
\label{lemma:test}
If the TFL is passed,
then the process is stable.
\end{proposition}

\subsection{Proof of Proposition~\ref{lemma:test}}
\label{sec:proof_test}

It is convenient to name the property used in the definition of the TFL.

\begin{property}
\label{prop:1}
We say that a function $ G=(G_i)_i:[0,\infty)\to \RR^K $ satisfies Property 
\ref{prop:1} for $\delta>0$ if \eqref{eq:verdura} holds for
every pair of times $ 0\le t_1<t_2 $ satisfying \eqref{eqn:ladrillo}.
\end{property}

Under this definition, the TFL is passed if and only if there exists $\delta>0$ such that every fluid limit satisfies Property \ref{prop:1} for $\delta>0$.

It is convenient to identify another similar property over functions.

\begin{property}
\label{prop:2}
We say that a function $ G=(G_i)_i:[0,\infty)\to \RR^K $ satisfies Property 
\ref{prop:2} for $\delta>0$ if \eqref{eq:verdura} holds for
every pair of times $ 0\le t_1<t_2 $ satisfying either
\begin{align}\label{eqn:cond1}
\min_{i\in L(t_1)}G_i(t)>\max_{i\notin L(t_1)}G_i(t)\quad \mbox{for every }t\in [t_1,t_2]
\end{align}
or
\begin{align}\label{eqn:cond2}
\min_{i\in L(t_2)}G_i(t)>\max_{i\notin L(t_2)}G_i(t)\quad \mbox{for every }t\in [t_1,t_2].
\end{align}
\end{property}

\begin{figure}[htbp]
\centering
\begin{minipage}{0.3\textwidth}
\centering
\resizebox{!}{3.3cm}{
	\begin{tikzpicture}
		\draw[thick,line cap=round] plot [smooth, tension=0.5] coordinates {(-2,1) (-1.5,0.7) (-1,1) (-0.5,0.5) (0.6,0.6) (1,0.1) (1.5,0.1) (2,0.5)};
		\draw[thick,line cap=round] plot [smooth, tension=0.5] coordinates {(-2,1) (-1.5,0.9) (-1,0.8) (-0.5,0.4)  (0.5,0.8) (0.8,-0.1) (1.5,0.2) (2,-0.5)};
		\draw[thick,line cap=round] plot [smooth, tension=0.5] coordinates {(-2,1) (-1.2,1.3) (-0.5,0.1)   (1,0.5)  (1.75,-0.1) (2,0)};
		\draw[gray,line width=0.5pt,line cap=round] plot [smooth, tension=0.5] coordinates {(-2,-0.5)  (-0.5,-0.5) (0.6,-1.5) (1.5,-1) (2,-1.7)};
		\draw[gray,line width=0.5pt,line cap=round] plot [smooth, tension=0.5] coordinates {(-2,-1.8) (-0.5,-1.4)  (0.5,-1.6)  (1.5,-1.7) (2,-1.5)};
		\draw[gray,line width=0.5pt,line cap=round] plot [smooth, tension=0.5] coordinates {(-2,-1) (-1.2,-0.7) (-0.5,-1)   (1,-1.5)  (2,-1)};
		\draw[gray,line width=0.5pt,line cap=round] plot [smooth, tension=0.5] coordinates {(-2,-1.5) (-1.2,-1.8)   (1,-1)  (2,-0.8)};
		\node[left] at (-2,1) {$L(t_1)$};
		\draw[line width=1pt] (-2,-2) -- (-2,-2.2);
		\node[below] at (-2,-2.2) {$t_1$};
		\draw[line width=1pt] (2,-2) -- (2,-2.2);
		\node[below] at (2,-2.2) {$t_2$};
	\end{tikzpicture}
}
\caption{Condition \eqref{eqn:cond1}}
\label{condition2a}
\end{minipage}
\hspace{30pt}
\begin{minipage}{0.3\textwidth}
\centering
\resizebox{!}{3.3cm}{
	\begin{tikzpicture}
		\draw[thick,line cap=round] plot [smooth, tension=0.5] coordinates {(-2,0.7) (-1.5,0.7) (-1,1.3) (-0.5,0.5) (0.6,0.6) (1,0.1) (1.5,0.1) (2,0)};
		\draw[thick,line cap=round] plot [smooth, tension=0.5] coordinates {(-2,1) (-1.5,0.9) (-1,0.8) (-0.5,0.4)  (0.5,0.8) (0.8,-0.1) (1.5,0.2) (2,0)};
		\draw[thick,line cap=round] plot [smooth, tension=0.5] coordinates {(-2,0.5) (-1.2,0.5) (-0.5,0.1)   (1,0.5)  (1.75,-0.1) (2,0)};
		\draw[gray,line width=0.5pt,line cap=round] plot [smooth, tension=0.5] coordinates {(-2,-0.5)  (-0.5,-0.5) (0.6,-1.5) (1.5,-1) (2,-1.7)};
		\draw[gray,line width=0.5pt,line cap=round] plot [smooth, tension=0.5] coordinates {(-2,-1.8) (-0.5,-1.4)  (0.5,-1.6)  (1.5,-1.7) (2,-1.5)};
		\draw[gray,line width=0.5pt,line cap=round] plot [smooth, tension=0.5] coordinates {(-2,-1) (-1.2,-0.7) (-0.5,-1)   (1,-1.5)  (2,-1)};
		\draw[gray,line width=0.5pt,line cap=round] plot [smooth, tension=0.5] coordinates {(-2,-1.5) (-1.2,-1.8)   (1,-1)  (2,-0.8)};
		\node[right] at (2,0) {$L(t_2)$};
		\draw[line width=1pt] (-2,-2) -- (-2,-2.2);
		\node[below] at (-2,-2.2) {$t_1$};
		\draw[line width=1pt] (2,-2) -- (2,-2.2);
		\node[below] at (2,-2.2) {$t_2$};
	\end{tikzpicture}
}
\caption{Condition \eqref{eqn:cond2}}
\label{condition2b}
\end{minipage}
\end{figure}

Conditions \eqref{eqn:cond1} and \eqref{eqn:cond2} are   represented in Figures~\ref{condition2a} and \ref{condition2b}, respectively.
Observe that condition~\eqref{eqn:ladrillo} in the definition of TFL implies both conditions \eqref{eqn:cond1} and \eqref{eqn:cond2}.
Said in different terms,
for a fixed $\delta>0$,
the range of pairs $t_1,t_2$ over which inequality \eqref{eq:verdura} is required to hold in Property~\ref{prop:2} is larger than in in Property~\ref{prop:1},
so every function satisfying Property~\ref{prop:2} also satisfies  Property~\ref{prop:1}.
The following  
result shows that the two properties indeed coincide over Lipschitz functions.

\begin{lemma}\label{lemma:equivalent_properties}
Let $ G:[0,\infty)\to\RR^k $ be a Lipschitz function,
and fix $\delta>0$.
Then $G$ satisfies Property \ref{prop:1} for $\delta$ if and only if it satisfies Property~\ref{prop:2} for $\delta$.
\end{lemma}

Before giving a proof of this lemma, we show how our methodological contribution, Proposition~\ref{lemma:test}, follows from it.
Suppose that the TFL is passed, that is, there exists $\delta>0$ such that every fluid limit satisfies Property \ref{prop:1} for this $\delta>0$.
Fix now a fluid limit $ G $.
By Lemma \ref{lemma:equivalent_properties}, $G$ satisfies Property \ref{prop:2} for $\delta$ because fluid limits are Lipschitz.
Fix $t$ such that $M(t)>0$ (recall that $M$ is the maximum function associated to $G$).
By continuity, there exists $ \varepsilon>0 $ such that \eqref{eqn:cond1} holds for $ t_1=t $ and any $ t_2\in (t,t+\varepsilon) $, which implies \eqref{eq:verdura}.
The analogous implication holds also to the left, i.e. for $ t_2=t $ and any $ t_1\in(t-\varepsilon,t) $, we have \eqref{eqn:cond2} and hence \eqref{eq:verdura}.
We have obtained
\begin{align}\label{eq:pizarra}
\limsup_{t'\to t}\frac{M(t)-M(t')}{t-t'}\le -\delta
\quad\mbox{for every $ t $ for which $ M(t)>0 $.}
\end{align}
Since $M$ is Lipschitz,
its derivative exists for almost every time,
and \eqref{eq:pizarra} proves that it is bounded by $-\delta$.
In other words, the maximum function is Lyapunov for the fluid limit, and the stability of the Markov process
under the SLC policy follows, see~\cite{Rob2003}.

\subsubsection*{Proof of Lemma \ref{lemma:equivalent_properties}}
We only prove the sufficiency direction since the necessity direction is trivial.
So assume that the Lipschitz function $G$ satisfies Property \ref{prop:1} for $\delta>0$.
We will proceed by induction over $m\in\NN$, the precise inductive hypothesis being the following one:
inequality \eqref{eq:verdura} holds
for $0\le t_1<t_2$
either  when $|L(t_1)|<m$ and \eqref{eqn:cond1} holds,
or when $|L(t_2)|<m$ and \eqref{eqn:cond2} holds.
Observe that the case $m=1$ holds because in this case Properties~\ref{prop:1} and \ref{prop:2} are equivalent.

Let now $0\le t_1<t_2$ such that $|L(t_1)|=m$ and \eqref{eqn:cond1} holds.
The trick is to define
\begin{align}
t^*=\max\{t\in [t_1,t_2]:L(t)=L(t_1)\}.
\end{align}
Since by assumption inequality~\eqref{eq:verdura} holds for every pair of times satisfying~\eqref{eqn:ladrillo}, we have  
\begin{align}\label{eqn:rosa}
\frac{M(t^*)-M(t_1)}{t^*-t_1}\le -\delta.
\end{align}
If $ t^*=t_2 $, we are trivially done.
Otherwise, for every $ t\in (t^*,t_2) $ we have that  $ |L(t)|< |L(t_1)| $ due to the validity of condition \eqref{eqn:cond1} and the very definition of $ t^* $.
Hence we can apply (both cases of) the inductive hypothesis to obtain that
\begin{align}
\limsup_{t'\to t}\frac{M(t')-M(t)}{t'-t}\le -\delta \quad\forall t\in (t^*,t_2).
\end{align}
From the fundamental theorem of calculus, we get
\begin{align}\label{eqn:jazmin}
\frac{M(t_2)-M(t^*)}{t_2-t^*}\le -\delta.
\end{align}
Inequalities \eqref{eqn:rosa} and \eqref{eqn:jazmin} let us obtain the desired inequality \eqref{eq:verdura}.
The second case in which $|L(t_2)|=m$ and \eqref{eqn:cond2}  follows analogously after defining
\begin{align}
t^*=\min\{t\in [t_1,t_2]:L(t)=L(t_2)\}.
\end{align}

\section{Proof of the maximum stability region, Theorem~\ref{thm:MSR}}
\label{sec:pmsr}


In this section we prove the sufficiency part of our main result,  Theorem \ref{thm:MSR}. 
More precisely, we will prove that the SLC policy is stable under conditions \eqref{eqMSRI}, \eqref{eq:condition_II} and \eqref{eqMSRIII} established for the corresponding settings.
To do so, in view of Proposition \ref{lemma:test}, we only need to verify that the corresponding fluid limits pass the TFL.
These verifications for Settings $\mathrm{II}$ and $\mathrm{III}$ are presented in Section~\ref{sec:TFLtwo} and Section~\ref{sec:TFLthree}, respectively,
while we relegate the verification for Setting~$\mathrm{I}$ to Appendix~\ref{subsection:pol_1}.
The proof  that the conditions \eqref{eqMSRI}, \eqref{eq:condition_II} and \eqref{eqMSRIII} are necessary for the existence of a stable policy is postponed to Appendix~\ref{app:nec}.

\label{subsection:pol_2_and_3}

We use the following standard queue length representation (see \cite{Dai1995} for instance).
Let $ A_i^x\subseteq [0,\infty) $ be the time subset in which queue $ i $ is receiving effective service.
By effective, we mean that we count only the time when the server is dedicated to this queue and the  queue is connected and non-empty.
It is formally defined as
\begin{align}\label{eqn:effective_time}
A_i^x=
\{t\in[0,\infty): C^x(t)=i,E_i^x(t)=1,Q_i^x(t)>0\}.
\end{align}
For every $ i \in[K]$, 
we call $\NNN_{\lambda_i}$ and $\NNN_{\mu_i}$ the counting measures  associated to the Poisson point processes 
respectively used to define the arrivals and departures in queue $i$.
The mentioned queue length representation is given by
\begin{align}\label{eq:representation_formula}
Q^x_i(t)=q_i+\NNN_{\lambda_i}([0,t])-\NNN_{\mu_i}(A^x_i\cap[0,t]),\quad i\in[K],t\ge 0.
\end{align}

Fix $G$ to be a fluid limit associated to an element $\omega\in\Omega'$ and a diverging sequence $(x^n)_n$.
Since $\omega$ is considered fixed, we do not include it in the subsequent notations.
An Arzelá-Ascoli argument allows us to prove a fluid version of formula \eqref{eq:representation_formula}, namely, for every $i\in[K]$, 
\begin{align}\label{eqn:lapicera}
G_i(t)=G_i(0)+\lambda_i t-\mu_i T_i(t).
\end{align}
Here $T_i$ is  the (uniformly over compacts) limit of the sequence of functions $(\bar T_i^{x^n})_n$ defined as
\begin{align}\label{eqn:holanda}
\bar T_i^{x^n}(t)=
\frac{\left|A^{x^n}_i\cap[0,t\|q^n\|]\right|}{\|q^n\|},
\quad t\ge 0.
\end{align}
This sequence represents the scaled versions of the cumulative effective service time.

\subsection{TFL  verification for Setting~II}
\label{sec:TFLtwo}
We first consider Setting~II. In order to prove that SLC is stable when~\eqref{eq:condition_II} holds, it is enough to show that the fluid limit corresponding to  the policy SLC satisfies the TFL. 
Recall the definitions of the max and argmax functions $M$ and $L$ associated to the fluid limit $G$.
Take $t_1$ and $t_2$ satisfying \eqref{eqn:ladrillo}.
For the sake of notational compactness, we call $L=L(t_1)$,
$\tilde\mu_i=\mu_i\pi_i(1)$ and $\tilde\rho_i=\lambda_i/\tilde\mu_i$.
We will prove that  \eqref{eq:verdura}, and hence TFL, is satisfied for
\begin{align}
\delta=\frac{1-\sum_{i=1}^K\tilde\rho_i}{\sum_{i=1}^K\tilde\mu_i^{-1}}.
\end{align}
Since this is a positive parameter in view of condition \eqref{eq:condition_II},
this implies the stability of the process through Proposition \ref{lemma:test}.
The numerator of $\delta$ has to be understood as the difference between the capacity of the server and the total effective load; the denominator is a term that appears when normalizing the queue sizes by the effective service rates.

Formula \eqref{eqn:lapicera}, in conjunction with assumption \eqref{eqn:ladrillo},  readily implies
\begin{align}
\frac{M(t_2)-M(t_1)}{t_2-t_1}\sum_{i\in L}
\tilde\mu_i^{-1}
=\sum_{i\in L}\frac{G_i(t_2)-G_i(t_1)}{\tilde\mu_i(t_2-t_1)}
=
\sum_{i\in L}\tilde\rho_i
-\sum_{i\in L}\frac{T_i(t_2)-T_i(t_1)}{\pi_i(1)(t_2-t_1)}.
\end{align}
We reduced the problem to proving that
\begin{align}\label{eqn:tenedor}
\sum_{i\in L}\frac{T_i(t_2)-T_i(t_1)}{\pi_i(1)(t_2-t_1)}=1
\end{align}
because, if true, we have
\begin{align}
\frac{M(t_2)-M(t_1)}{t_2-t_1}=
\frac{\sum_{i\in L}\tilde\rho_i-1}{\sum_{i\in L}
\tilde\mu_i^{-1}}\le -\delta.
\end{align}

We explain now why identity \eqref{eqn:tenedor} holds,
and we refer for the details to Appendix~\ref{sct:calor}.
The first observation is that, since we are assuming \eqref{eqn:ladrillo},
the queues in $L$ will have the leading queue sizes during all the considered time interval in the pre-limit.
Also, the non-preemptive nature of the policy guarantees that there is always a connected queue at every decision epoch.
These two observations imply that, in this regime in which the queues in $L$ are leading, the server will be dedicated precisely to this queue subset~$L$ (except possibly during an initial negligible time interval). 
The LHS of Equation~(\ref{eqn:tenedor}) captures the proportion of time that the server is assigned to queues in $L$ during the time interval $[t_1,t_2]$. To see this, we note that $T_i(t_2)-T_i(t_1)$ is the amount of time the server is dedicated to $i$ while this queue is connected. When we divide this time by the proportion  of time in which queue $i$ is effectively connected during this interval, which is precisely $\pi_i(1)$ as explained in the next paragraph, we get the total amount of time that the server is dedicated to queue~$i$. Thus, $ \frac{T_i(t_2)-T_i(t_1)}{\pi_i(1)(t_2-t_1)}$ represents the proportion of time that the server is dedicated to $i$. 
Since, as just explained, the server is fully dedicated to  queues in $L$, we can conclude that these fractions sum up to 1, yielding (\ref{eqn:tenedor}).


We now explain why the proportion described above is given by $\pi_i(1)$.
Consider the stochastic process that encodes the state of the $i$-th environment only during the time the server is dedicated to it.
This is a Markov process with state-space $\{0,1\}$ that transitions from disconnected to connected  at rate $\lambda_i'$, and from connected to disconnected at  rate $\mu_i'$.
Its invariant distribution is precisely $\pi_i$, and by the ergodic theorem $\pi_i(1)$ is the desired proportion.
For clarity, it is important to note that,
if a departure  occurs in queue~$i$,  
the state of its environment, 1,  remains unchanged.
Indeed,
since the policy is in the family $\PPP$,
even if the server is dedicated to other queues for a period of time after this departure, the $i$-th queue will necessarily be connected again when the  server is re-dedicated to it.

\subsection{TFL  verification for Setting~III}
\label{sec:TFLthree}
Define $\hat\mu_i=\theta_i(\gamma)\mu_i$ and $\hat\rho_i=\lambda_i/\hat\mu_i$ for notational compactness.
As in the previous case, but replacing $\pi_i(1)$ by $\theta_i(\gamma)$,
we get
\begin{align}
\frac{M(t_2)-M(t_1)}{t_2-t_1}\sum_{i\in L}
\hat\mu_i^{-1}
=
\sum_{i\in L}\hat\rho_i
-\sum_{i\in L}\frac{T_i(t_2)-T_i(t_1)}{\theta_i(\gamma)(t_2-t_1)}.
\end{align}
We will prove in Appendix~\ref{sct:calor} that
\begin{align}\label{eqn:prisma}
\sum_{i\in L}\frac{T_i(t_2)-T_i(t_1)}{\theta_i(\gamma)(t_2-t_1)}=1-\pi_L^0,
\end{align}
which implies the desired inequality \eqref{eq:verdura} for
\begin{align}
\delta=\frac{\min\{1-\pi_L^{0}-\sum_{i\in L}\hat\rho_i:L\subseteq[K],L \neq\emptyset\}}
{\sum_{i=1}^K\hat\mu_i^{-1}}.
\end{align}
The numerator of $\delta$ represents the minimum over subsets~$L$ of the difference between the maximum capacity that the server can dedicate to the subset and the total effective load in question.

Identity \eqref{eqn:prisma} has a similar explanation as \eqref{eqn:tenedor}, with some differences that we now outline. For further details we refer to Appendix~\ref{app:g}. 
The first difference is that the coefficient equals $\theta_i(\gamma)$ and not $\pi_i(1)$.
As before, we define the Markov process encoding the state of the $i$-th environment during the time the server is dedicated to queue~$i$,
which in this case transitions from disconnected to connected at rate $\lambda'_i+\gamma$, and from connected to disconnected at rate $\mu_i'$.
The reason is akin to what was mentioned earlier, with the additional matter that we pass from disconnected to connected also when the exponential clock of intensity $\gamma$ rings.
Hence, by the ergodic theorem,
$\theta_i(\gamma)=\frac{\lambda_i'+\gamma}{\mu_i'+\lambda_i'+\gamma}$
is the asymptotic proportion of time this process is connected as expected.

The second difference is that the RHS of \eqref{eqn:prisma} is $1-\pi_L^0$ and not $1$ as in the previous case.
To understand this consider two consecutive $\gamma$-marks $s_1,s_2\in (t_1^n,t_2^n]$.
In an interval $[s_1,s_2)$ the server is either dedicated to a fixed queue, or is dedicated to $\emptyset$. Since we are under Setting~III and consider SLC,  a queue~$i$ in $L$ is always served, except when all queues in $L$ are disconnected at time~$s_1$. The latter happens with probability $\pi_L^0$. 
Therefore, 
$1-\pi_L^0$ is the asymptotic proportion of this type of intervals (those within $(t_1^n,t_2^n]$ that are delimited by $\gamma$-marks) for which queues in~$L$ are  served. Since the LHS of~\eqref{eqn:prisma} is the proportion of time that the server is assigned to queues in~$L$, this should hence be equal to $1-\pi_L^0$. 

\section{Numerical evaluation of stability regions}
\label{sec:num}

In this section we provide some numerical examples of the comparison results we described in Section~\ref{sec:main}. In particular we want to assess the impact of $\gamma$ on the reduction in the stability region, $\SR{III}(\gamma)$, for which a closed-form expression is given in~\eqref{eq:redIII}.

\begin{figure}[t]
\begin{center}
\includegraphics[width=0.48\textwidth]{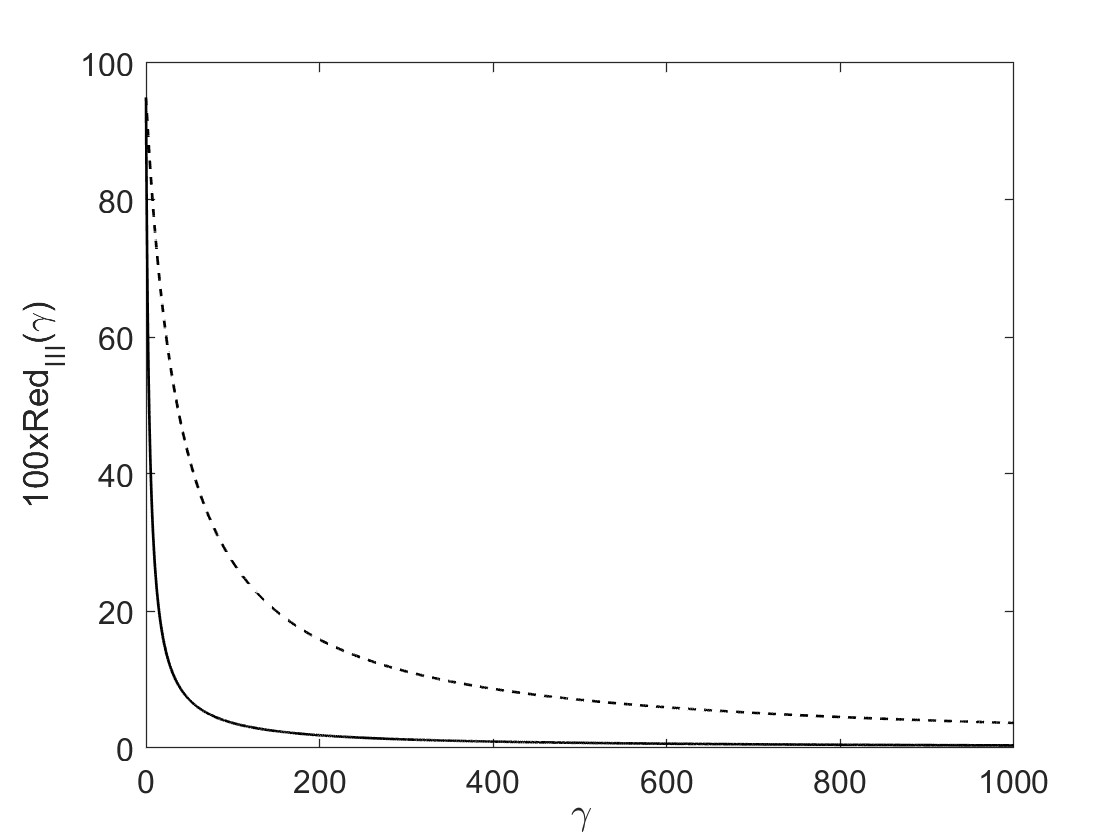}
\includegraphics[width=0.48\textwidth]{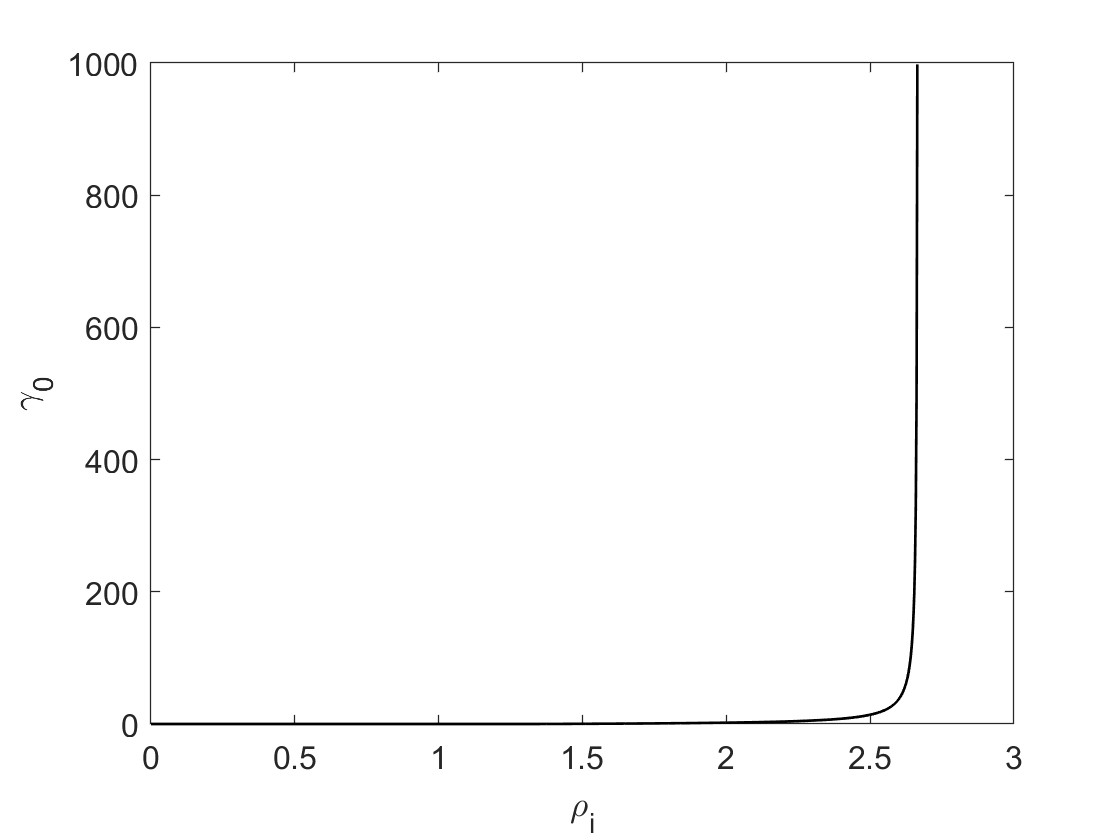}
\caption{ {\em(left)}  $\SR{III}(\gamma)$ as a function of the parameter $\gamma$, with $\lambda_1'=\alpha 0.8, \lambda_2'=\alpha 3, \mu_1'=0.2\alpha $ and $\mu_2'=\alpha$. We set $\alpha=1$ for the \emph{solid} line and $\alpha=10$ for the \emph{dashed line}. 
	{\em (right)} Minimal $\gamma_0$ as a function of $\rho_i$, with $\lambda_i'=\mu_i'=1$ and $\rho_1=\rho_2$.  }
\label{fig:reducvsgamma}
\end{center}
\end{figure}

In Figure~\ref{fig:reducvsgamma}{ \em (left)}, we plot the reduction in the maximum stability region when decision are at $\gamma$ epochs, i.e.\ $\SR{III}(\gamma)$, as a function of the decision rate $\gamma$. 
The difference between the two lines lies in the speed of the environments, which is 10 times faster in the \emph{dash} line.
As expected from~\eqref{eq:redIII}, the reduction of the stability decreases as $\gamma$ increases and converges to zero. We further observe that the \emph{solid} line decreases much faster. The intuition for this is that the state of connectivity of the queues changes less fast when $\alpha=1$, and hence, a smaller rate of decisions is already sufficient to get  a similar stability reduction as for the case where the connectivity changes faster, $\alpha=10$.

In Figure~\ref{fig:reducvsgamma} {\em (right)} we illustrate the result of Corollary~\ref{cor:gamma0}. We take a  symmetric situation with $\lambda_i'=\mu_i'=1$ and  $\rho_1=\rho_2$. From Theorem~\ref{thm:MSR} we obtain that in order to be stable in Setting~I, the load should satisfy $\rho_i \in [0,(1-\pi_1(0)\pi_2(0))/2) $.   In Figure~\ref{fig:reducvsgamma} {\em (right)} we plot $\gamma_0$ as a function of $\rho_i$. We recall that $\gamma_0$, defined in Corollary~\ref{cor:gamma0}, denotes the minimal value of $\gamma$ such that the system is stable in Setting~III. We  observe that the minimal value of $\gamma_0$ remains relatively small until $\rho_i$ gets very close to the stability border of $\MSR{I}$, which shows  a sharp phase transition: stabilizing the system in terms of decision rate is relatively cheap until we are very close to the frontier.  

\begin{figure}[t]
\begin{center}
\includegraphics[width=0.48\textwidth]{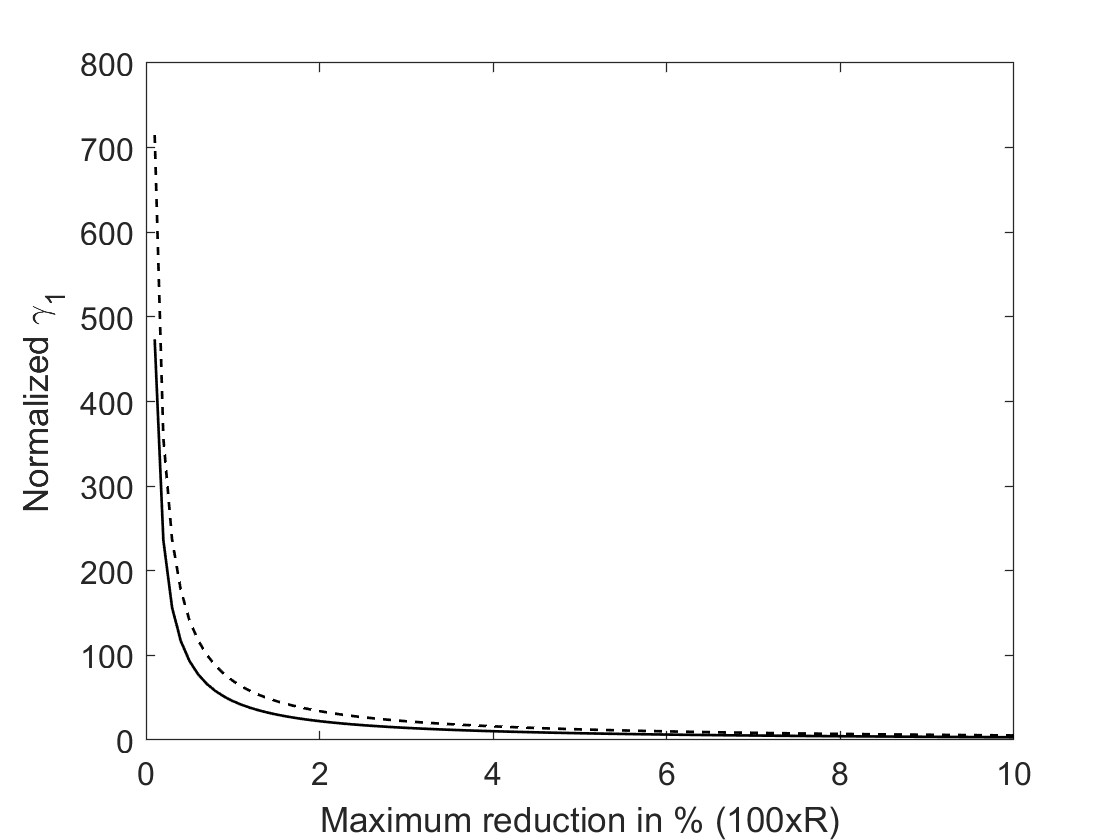}
\includegraphics[width=0.48\textwidth]{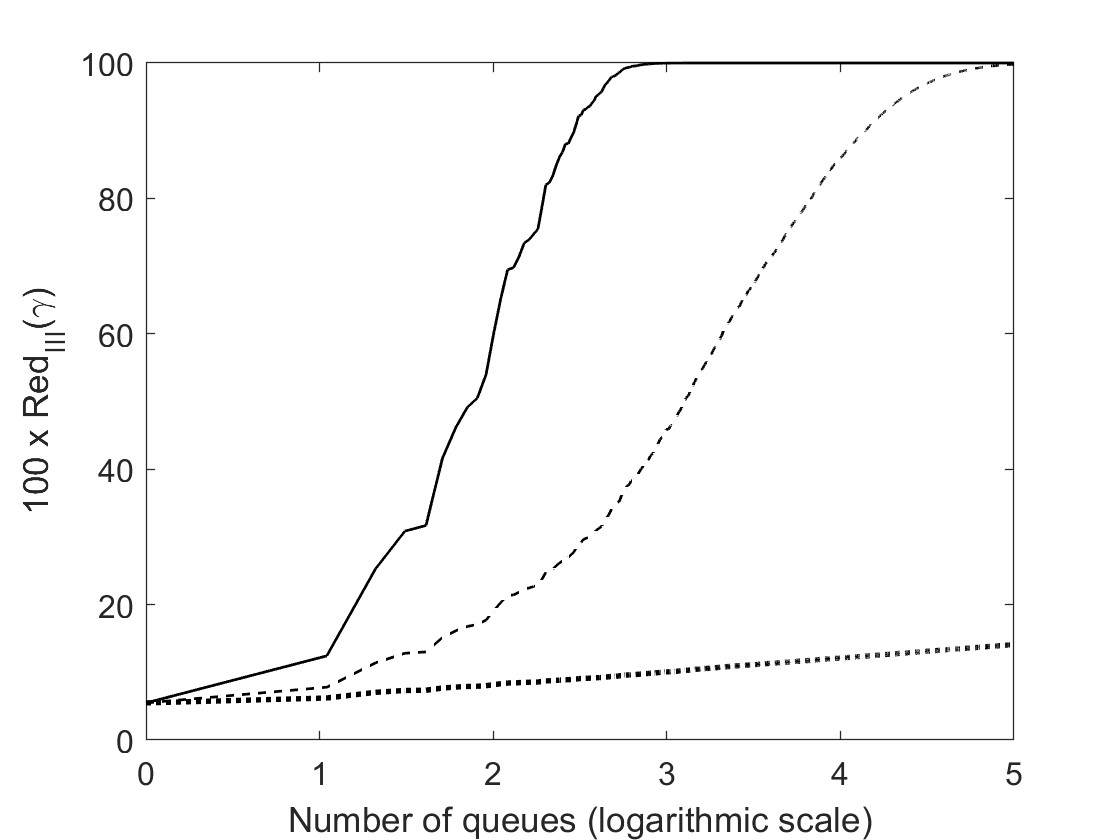}
\caption{
{\em (left)}
The normalized value  $\frac{\gamma_1-\mbox{U}}{\mbox{U}}$ as a function of the tolerable stability reduction $R$. 
The parameters from Figure~\ref{fig:reducvsgamma} \emph{(left)} are used here.
{\em (right)} $\SR{III}(\gamma)$ versus the number of queues, where $\gamma= 5$  for the \emph{solid} line, $\gamma = 5\sqrt{K}$  for the \emph{dashed} curve and $\gamma= 5 K$ for the \emph{dotted} line. 
} 
\label{fig:reducvsservers}
\end{center}
\end{figure}

\begin{figure}[t]
\begin{center}
\includegraphics[width=0.48\textwidth]{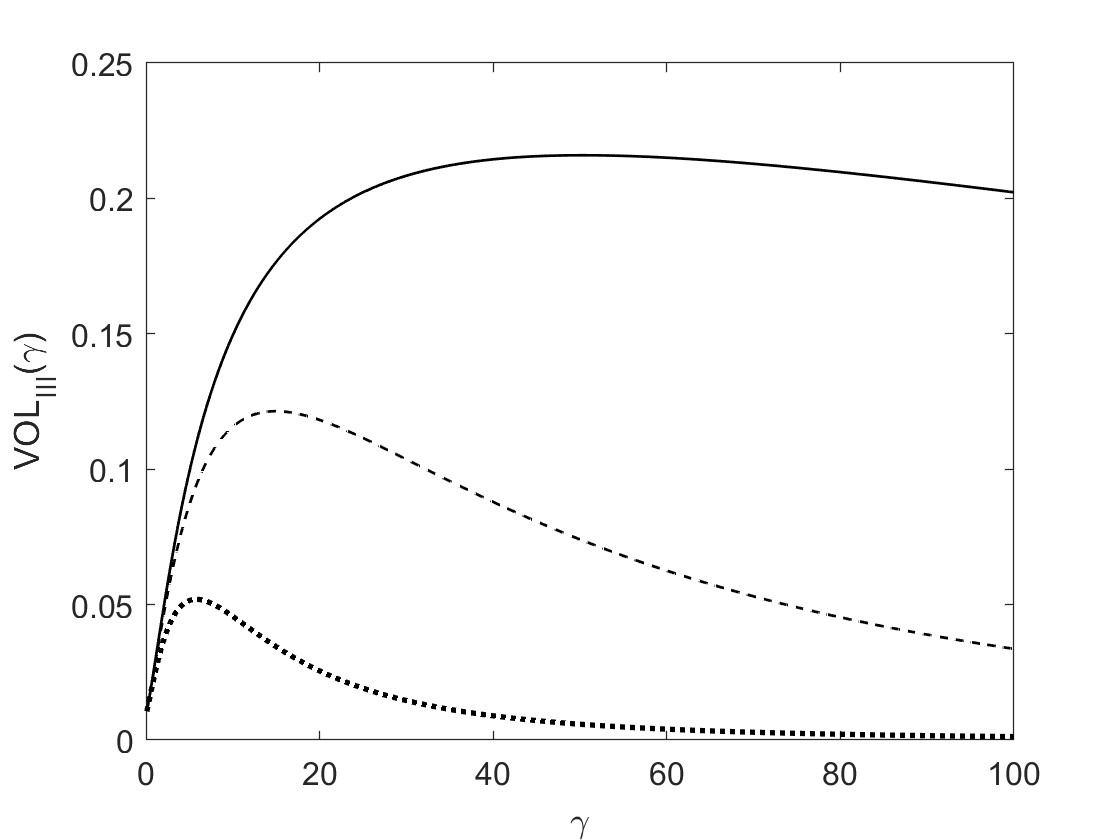}
\includegraphics[width=0.48\textwidth]{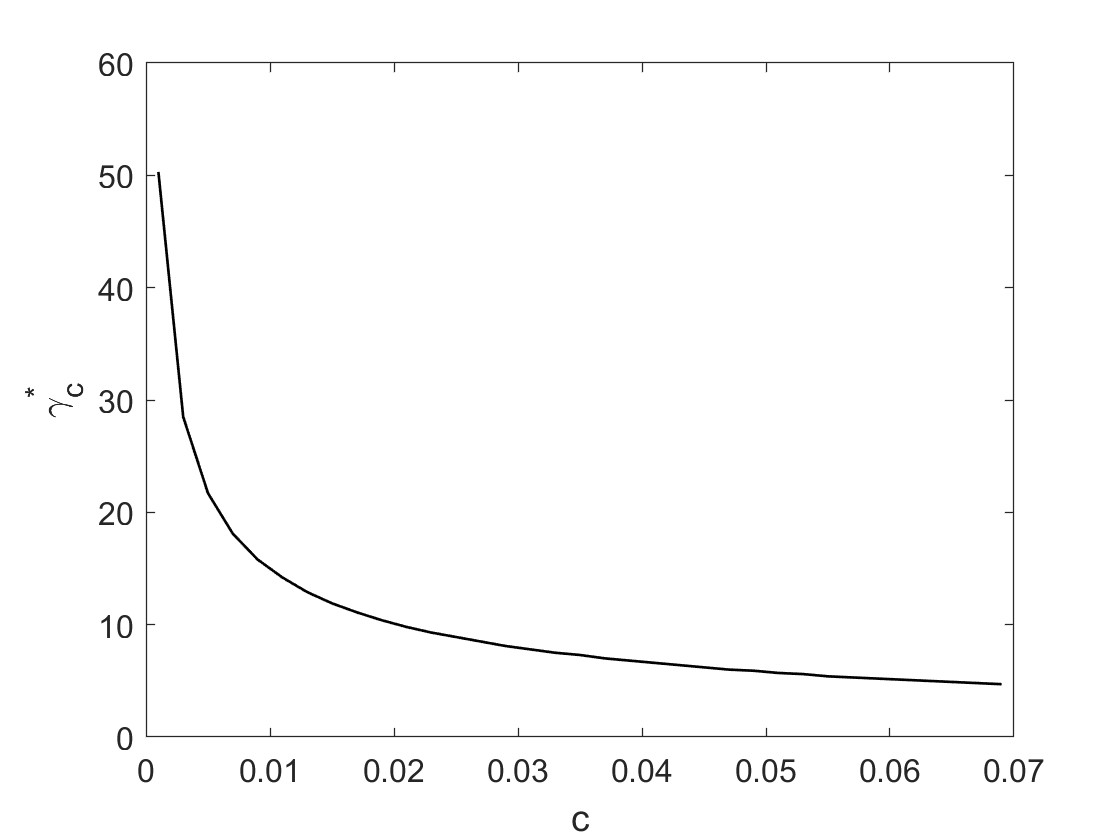}
\caption{ We set    $\lambda_1'=1$, $\lambda_2'=2$, $\lambda_3'=3$, $\mu_1'=3$ and $\mu_2'=2$ and $\mu_1'=1$. 
{\em (left)}
$\VOL{III}(\gamma)$ (as determined by the inequalities~(\ref{eq:optrateoverheads})) as a function of $\gamma$. The inactivity duration $c$ equals $c=0.001$ for the \emph{solid} line,  $c=0.01$ for the \emph{dashed} line, and $c=0.05$ for the \emph{dotted} line.
{\em (right)} Optimal decision rate, $\gamma_c^*$, as a function of $c$.
} 
\label{fig:optrate}
\end{center}
\end{figure}

In Figure~\ref{fig:reducvsservers} {\em (left)} we  illustrate the result of Corollary~\ref{cor:gam1}, that is, we plot the minimum rate $\gamma_1$ such that for any $\gamma\geq \gamma_1$, the reduction in stability for Setting~III is less than the tolerable reduction given by~$R$. 
To find $\gamma_1$, we simply solve the 2nd degree polynomial obtained from $\SR{III}(\gamma) = (1 -\theta_1(\gamma) \theta_2(\gamma))=R$. 
Again, we consider two set of parameters, where the difference is that the rates of connectivity are scaled. In order to compare the values of $\gamma_1$, we scale them by $\mbox{U} = \lambda_1'+\lambda_2'+\mu_1'+\mu_1'$ and  plot  $(\gamma_1-\mbox{U})/\mbox{U}$. 
Firstly, we see that both curves overlap quite closely, despite the fact that the rates of the environments for the dashed line are 10 times larger. 
Secondly, we observe that if we want to make the stability reduction to be smaller than, say 1$\%$,  the parameter $\gamma_1$ needs to be around 50 times larger than the sum of the transition rates, $U$.

In Figure~\ref{fig:reducvsservers} {\em(right)} we plot the reduction in stability as a function of the number of queues, $K$. The parameters of the  queue,  $\lambda_i'$ and $\mu_i'$, are chosen randomly from a uniform distribution on $[0,1]$. 
In the \emph{solid} curve, the value of $\gamma$ is kept constant, and we see that the reduction tends to 100$\%$, as could also be seen directly from Equation~(\ref{eq:redIII}).
In the \emph{dash} line, we scale the value of $\gamma$ with the square root of the number of servers, and we observe that, even though more slowly, the reduction still converges to 100$\%$. In the \emph{dotted} line, we scale $\gamma$ linearly with the number of servers and we observe that the stability reduction tends to some positive value strictly smaller than 1. 
It would be interesting to show that indeed scaling the decision rate $\gamma$ linearly with the number of servers yields a constant stability reduction asymptotically.

In Figure~\ref{fig:optrate}, we illustrate the result of Corollary~\ref{optigamma} where the stability region is determined for a system where each decision epoch causes an inactivity period of duration~$c$. We set $K=3$ and in Figure~\ref{fig:optrate} {\em (left)} we plot $\VOL{III}(\gamma)$ (can be calculated from   Equation~\ref{eq:optrateoverheads}) as a function of~$\gamma$. We observe the existence of an optimal rate $\gamma_c^*$ that maximizes the stability region and hence strikes the right balance between a high frequency of taking actions and the overheads induced by them. 
In Figure~\ref{fig:optrate}~{\em (right)} we plot $\gamma_c^*$ as a function of~$c$. We observe that $\gamma_c^*$ is decreasing in~$c$. Indeed, the larger $c$, the more harmful it becomes to increase the frequency of actions.

\section{Conclusions}

Building on an original analysis of fluid limits of one server multi-class system in a random environment,
our results reveal the crucial impact of decision epochs on the stability properties of a natural class $\PPP$ of policies. We conclude mentioning a few research avenues that stem out of our work.
Exploring the maximal stability region for policies
outside $\PPP$ remains an open problem (see Section \ref{sec:SLC_is_not_MS}).
Besides, exploring the impact of decision epochs on the stability region opens avenues for further investigation, for instance in stochastic networks with more complex topologies and broader statistical assumptions.

More generally, our approach finds potential application for the study of continuous-time Partially Observable Markov Decision Processes (POMDP) with constrained action epochs and their model-free counterpart in reinforcement learning, which currently gathers considerable attention. 
Finally, we also mention that the methodological contribution we present on fluid limits, could be promising for examining the stability of various other models.


\bibliographystyle{plain}
\bibliography{biblio}

\begin{thebibliography}{10}

\bibitem{andrews2004scheduling}
M.~Andrews, K.~Kumaran, K.~Ramanan, A.~Stolyar, R.~Vijayakumar, and P.~Whiting.
\newblock Scheduling in a queuing system with asynchronously varying service
  rates.
\newblock {\em Probability in the Engineering and Informational Sciences},
  18(2):191--217, 2004.

\bibitem{ayesta2013}
U.~Ayesta, M.~Erausquin, M.~Jonckheere, and I.~M. Verloop.
\newblock Scheduling in a random environment: Stability and asymptotic
  optimality.
\newblock {\em IEEE/ACM Transactions on Networking}, 21(1):258--271, 2013.

\bibitem{BacceliMakowski}
F.~Baccelli and A.~M. Makowski.
\newblock Stability and bounds for single server queues in random environment.
\newblock {\em Communications in Statistics. Stochastic Models}, 2(2):281--291,
  1986.

\bibitem{bambos2004queueing}
N.~Bambos and G.~Michailidis.
\newblock Queueing and scheduling in random environments.
\newblock {\em Advances in Applied Probability}, 36(1):293--317, 2004.

\bibitem{BM2005}
Nicholas Bambos and George Michailidis.
\newblock Queueing networks of random link topology: stationary dynamics of
  maximal throughput schedules.
\newblock {\em Queueing Syst.}, 50(1):5--52, 2005.

\bibitem{Borst05}
S.C. Borst.
\newblock User-level performance of channel-aware scheduling algorithms in
  wireless data networks.
\newblock {\em IEEE/ACM Transactions on Networking}, 13(3):636--647, 2005.

\bibitem{Bra2008}
M.~Bramson.
\newblock {\em Stability of queueing networks}, volume 1950.
\newblock Springer, Berlin, Berlin, 2008.

\bibitem{budhiraja14}
Amarjit Budhiraja, Arka Ghosh, and Xin Liu.
\newblock Scheduling control for markov-modulated single-server multiclass
  queueing systems in heavy traffic.
\newblock {\em Queueing Systems}, 78(1):57--97, 2014.

\bibitem{celik2011scheduling}
G.~D Celik, L.~B. Le, and E.~Modiano.
\newblock Scheduling in parallel queues with randomly varying connectivity and
  switchover delay.
\newblock In {\em 2011 Proceedings IEEE INFOCOM}, pages 316--320. IEEE, 2011.

\bibitem{Dai1995}
J.~G. Dai.
\newblock On positive {H}arris recurrence of multiclass queueing networks: a
  unified approach via fluid limit models.
\newblock {\em Ann. Appl. Probab.}, 5(1):49--77, 1995.

\bibitem{Dauria}
B.~D'Auria.
\newblock $m/m/\infty$ queues in semi-markovian random environment.
\newblock {\em Queueing Systems}, (58):221--237, 2008.

\bibitem{dimakis_walrand_2006}
Antonis Dimakis and Jean Walrand.
\newblock Sufficient conditions for stability of longest-queue-first
  scheduling: second-order properties using fluid limits.
\newblock {\em Advances in Applied Probability}, 38(2):505–521, 2006.

\bibitem{dinh2013survey}
H.~T Dinh, C.~Lee, D.~Niyato, and P.~Wang.
\newblock A survey of mobile cloud computing: architecture, applications, and
  approaches.
\newblock {\em Wireless communications and mobile computing},
  13(18):1587--1611, 2013.

\bibitem{DURANSIGMETRICS}
Santiago Duran and Ina Maria~Maaike Verloop.
\newblock {Asymptotic Optimal Control of Markov-Modulated Restless Bandits}.
\newblock In {\em {ACM Sigmetrics 2018}}, Irvine, United States, June 2018.

\bibitem{khabbaz2015modelling}
M.~Khabbaz and C.~M. Assi.
\newblock Modelling and analysis of a novel deadline-aware scheduling scheme
  for cloud computing data centers.
\newblock {\em IEEE Transactions on Cloud Computing}, 6(1):141--155, 2015.

\bibitem{Kleinrock76_1}
L.~Kleinrock.
\newblock {\em Queueing Systems, vol. 1}.
\newblock John Wiley and Sons, 1976.

\bibitem{ModianoLQF}
L.~B. Le, E.~Modiano, C.~Joo, and N.~B. Shroff.
\newblock Longest-queue-first scheduling under the sinr interference model.
\newblock In {\em ACM MobiHoc}, 2010.

\bibitem{peng2015random}
Z.~Peng, D.~Cui, J.~Zuo, Q.~Li, B.~Xu, and W.~Lin.
\newblock Random task scheduling scheme based on reinforcement learning in
  cloud computing.
\newblock {\em Cluster computing}, 18:1595--1607, 2015.

\bibitem{Rob2003}
P.~Robert.
\newblock {\em Stochastic networks and queues}, volume~52 of {\em Applications
  of Mathematics (New York)}.
\newblock Springer-Verlag, Berlin, french edition, 2003.
\newblock Stochastic Modelling and Applied Probability.

\bibitem{shakkottai2002scheduling}
S.~Shakkottai and A.~L. Stolyar.
\newblock Scheduling for multiple flows sharing a time-varying channel: The
  exponential rule.
\newblock {\em Translations of the American Mathematical Society-Series 2},
  207:185--202, 2002.

\bibitem{srikant2013communication}
R.~Srikant and L.~Ying.
\newblock {\em Communication Networks: An Optimization, Control, and Stochastic
  Networks Perspective}.
\newblock Cambridge University Press, 2013.

\bibitem{Stolyar04}
A.L. Stolyar.
\newblock Maxweight scheduling in a generalized switch: State space collapse
  and workload minimization in heavy traffic.
\newblock {\em Annals of Applied Probability}, 14(1):1--53, 2004.

\bibitem{tassiulas1993dynamic}
L.~Tassiulas and A.~Ephremides.
\newblock Dynamic server allocation to parallel queues with randomly varying
  connectivity.
\newblock {\em IEEE Transactions on Information Theory}, 39(2):466--478, 1993.

\bibitem{Verloop16}
I.M. Verloop.
\newblock Asymptotically optimal priority policies for indexable and
  non-indexable restless bandits.
\newblock {\em Annals of Applied Probability}, 26(4):1947--1995, 2016.

\bibitem{weber1990index}
Richard~R Weber and Gideon Weiss.
\newblock On an index policy for restless bandits.
\newblock {\em Journal of applied probability}, pages 637--648, 1990.

\bibitem{Whi88}
P.~Whittle.
\newblock Restless bandits: Activity allocation in a changing world.
\newblock {\em Journal of Applied Probability}, 25:287--298, 1988.

\end{thebibliography}

\newpage
\appendix

\section{Necessity proof of Theorem \ref{thm:MSR}}
\label{app:nec}

We prove now that the conditions given in Theorem \ref{thm:MSR} are necessary.

\subsubsection*{Setting I}
For this setting, this result is not new (see \cite{BM2005} for instance),
but nevertheless we give a proof for completeness.
Suppose that condition \eqref{eqMSRI} is not fulfilled,
namely,
that there exists a non-empty subset of queues $L\subseteq[K]$ such that
\begin{align}
\sum_{i\in L}\rho_i>1-\pi_{L}^0.
\end{align}
Under this condition, we will prove that any policy is unstable,
the idea being that, even in the extreme scenario in which all the service is dedicated to the subset $L$,
tasks will necessarily accumulate in these queues.

From formula \eqref{eq:representation_formula},
for an arbitrary initial condition $x$,
we get
\begin{align}\label{eq:andes}
\sum_{i\in L}\frac{Q^x_i(t)}{t\mu_i}=
\sum_{i\in L}\frac{q_i}{t\mu_i}
+\sum_{i\in L}\frac{\NNN_{\lambda_i}([0,t])}{t\mu_i}
-\sum_{i\in L}\frac{\NNN_{\mu_i}(A^x_i\cap[0,t])}{t\mu_i}.
\end{align}
The last sum can be controlled as follows:
\begin{align}
\sum_{i\in L}\frac{\NNN_{\mu_i}(A^x_i\cap[0,t])}{t\mu_i}
&\le\sum_{i\in L}\frac{\NNN_{\mu_i}\left(\left\{s\in[0,t]:C^x(s)=i,E^x_i(s)=1\right\}\right)}{t\mu_i}
\\ \label{eqn:kiwi}
&=\frac{1}{t}\sum_{i\in L}\left|\left\{s\in[0,t]:C^x(s)=i,E^x_i(s)=1\right\}\right|+R_t
\\
&\le \frac{1}{t}\Big|\bigcup_{i\in L}\left\{s\in[0,t]:E^x_i(s)=1\right\}\!\Big|+R_t.
\end{align}
In the first inequality, we simply neglected the per-queue unbusy periods, i.e., the periods in which the queue is empty.
The error $R_t$ is defined for the equality to hold.
The second inequality involves two steps:
first,
we use that the sets appearing in the sum in \eqref{eqn:kiwi} are disjoint;
secondly, as anticipated, we bound by the most beneficial case in which all the service is dedicated to the queues in $L$.
Since $R_t$ goes to zero as $t\to\infty$ ---due to the law of large numbers---,
the last expression converges almost surely to $1-\pi_L^0$.
Coming back to \eqref{eq:andes}, and using that the second sum in the RHS of \eqref{eq:andes} converges to $\sum_{i\in L}\rho_i$,
we obtain
\begin{align}
\liminf_{t\to\infty}\sum_{i\in L}\frac{Q^x_i(t)}{t\mu_i}
\ge \sum_{i\in L}\rho_i-(1-\pi_L^0)
\quad a.s.
\end{align}
Since by assumption the RHS in this inequality is strictly positive, we can conclude that
\begin{align}
\liminf_{t\to\infty}
\sum_{i\in L}\frac{Q_i(t)}{\mu_i}=+\infty\quad a.s.,
\end{align}
and hence the process is unstable.

\subsubsection*{Setting II}
Suppose that condition \eqref{eq:condition_II} is not satisfied,
i.e., suppose that
\begin{align}
\sum_{i=1}^K\frac{\rho_i}{\pi_i(i)}>1.
\end{align}
Under this assumption, we will prove that any policy is unstable.

Again from formula \eqref{eq:representation_formula},
\begin{align}\label{eq:pirineos}
\sum_{i=1}^K\frac{Q_i(t)}{t\pi_i(1)\mu_i}=
\sum_{i=1}^K\frac{q_i}{t\pi_i(1)\mu_i}
+\sum_{i=1}^K\frac{\NNN_{\lambda_i}([0,t])}{t\pi_i(1)\mu_i}
-\sum_{i=1}^K\frac{\NNN_{\mu_i}(A^x_i\cap[0,t])}{t\pi_i(1)\mu_i}.
\end{align}
We can deal with the last sum as follows:
\begin{align}
\sum_{i=1}^K\frac{\NNN_{\mu_i}(A^x_i\cap[0,t])}{t\pi_i(1)\mu_i}
&\lessapprox
\frac{1}{t}\sum_{i=1}^K
\frac{\left|\left\{s\in[0,t]:C^x(s)=i,E^x_i(s)=1\right\}\right|}{\pi_i(1)}
\\ \label{eqn:conejo}
&\lessapprox
\frac{1}{t}\sum_{i=1}^K
\left|\left\{s\in[0,t]:C^x(s)=i\right\}\right|.
\end{align}
The symbol $\lessapprox$ means `less or equal  with an error that vanishes as $t\to\infty$'.
In the first line, we neglected the unbusy periods and used the law of large numbers, 
exactly as we did right before in the necessity proof for Setting I.
To justify the second line,
we refer to the argument that we used at the end of the test verification for the same Setting II, in Section \ref{subsection:pol_2_and_3}.
In that argument,
we used that $\pi_i(1)$ is the asymptotic proportion of time in which a certain queue $i$ is connected within the time it is receiving service.
In the case of a general policy, $\pi_i(1)$ instead works as an upper bound for such a proportion because it may happen that, at a decision epoch, we serve a queue that is disconnected.
If this happens, the time interval between such epoch and the moment in which the queue under service connects, is a disconnected time interval that is not counted in the SLC case.
Using that the expression in \eqref{eqn:conejo} is bounded by $1$,
we obtain
\begin{align}
\limsup_{t\to\infty}\sum_{i=1}^K\frac{\NNN_{\mu_i}(A^x_i\cap[0,t])}{t\pi_i(1)\mu_i}\le 1\quad a.s.
\end{align}
Substituting in \eqref{eq:pirineos}, we get
\begin{align}
\liminf_{t\to\infty}
\sum_{i=1}^K\frac{Q_i(t)}{t\mu_i\pi_i(1)}
\ge
\sum_{i=1}^K\frac{\rho_i}{\pi_i(i)}-1
\quad a.s.,
\end{align}
which let us conclude.

\subsubsection*{Setting III}
Suppose that 
\begin{align}\label{eqn:dado}
\sum_{i\in L}\frac{\rho_i}{\theta_i(\gamma)}>1-\pi_{L}^0
\end{align}
for some non-empty subset $L\subseteq[K]$.
Unlike the other settings,
in Setting III we can only prove that any policy in the family $\PPP$ is unstable.
This is related with the discussion given in \ref{sec:SLC_is_not_MS}, in which me mention that, even if \eqref{eqn:dado} holds, there might stable policies outside the family $\PPP$.

The proof in this setting has the same spirit than in the previous ones.
In this case we have again formula \eqref{eq:pirineos},
but with $\theta_i(\gamma)$ instead of $\pi_i(1)$,
and with the sums running only over the queues in $L$.
We control the sum concerning the services as
\begin{align}
\sum_{i\in L}\frac{\NNN_{\mu_i}(A^x_i\cap[0,t])}{t\theta_i(\gamma)\mu_i}
&\lessapprox
\frac{1}{t}\sum_{i\in L}
\frac{\left|\left\{s\in[0,t]:C^x(s)=i,E^x_i(s)=1\right\}\right|}{\theta_i(\gamma)}
\\ \label{eqn:conejo}
&\approx
\frac{1}{t}\sum_{i\in L}
\left|\left\{s\in[0,t]:C^x(s)=i\right\}\right|
\\
&\lessapprox
1-\pi_L^0.
\end{align}
Of course, the symbol $\approx$ means `equal  with an error that vanishes as $t\to\infty$'.
In the first line, we neglected the per-queue unbusy periods, as usual.
In the second line, we used that $\theta_i(\gamma)$ is the proportion of connected periods within the time that we are serving queue $i$.
We emphasize that this is an inherent property of the family $\PPP$ and not of the SLC policy, and this is also the reason why, in this line, we have an $\approx$ instead of an $\lessapprox$ as we had in the Setting II.
The last line is due to the fact that the proportion of time intervals separated by consecutive $\gamma$-marks in which we serve a queue in $L$
is at most $1-\pi_L^0$  because, since the policy is in the family $\PPP$, at decision epochs we cannot select a queue that is disconnected.
Importantly, this last step would fail if we considered policies outside $\PPP$.

\section{Remaining proofs from the TFL verifications}

\subsection{Proof of Equations (\ref{eqn:tenedor})}
\label{sct:calor}

This identity  follows from the convergence
\begin{align}\label{eqn:madrigal}
\lim_{n\to\infty}\sum_{i\in L}\frac{\bar T_i^{x^n}(t_2) -\bar T_i^{x^n}(t_1)}{\pi_i(1)(t_2-t_1)}=1
\end{align}
simply because, by definition, $\bar T_i^{x^n}$ approximates $T_i$.
Let $t_1^n=\|q^n\|t_1$ and $t_2^n=\|q^n\|t_2$ represent the microscopic versions of $t_1$ and $t_2$.
Let also $\tau_n$ be the first time after $t_1^n$ a queue in $L$ receives service, formally defined as
\begin{align}\label{eqn:enchufla}
\tau_n=\inf\{t\ge t_1^n:C^{x^n}(t)\in L\}.
\end{align}
Considering that it is sufficient for one of the queues in $L$ to be connected in an decision epoch for the SLC policy to start serving this group of queues, it is not difficult to believe that
\begin{align}\label{eqn:cebolla}
\lim_{n\to\infty}\frac{\tau_n}{\|q^n\|}=t_1.
\end{align}
This is proven at the end of the section.

Convergence \eqref{eqn:madrigal} follows from the following sequence of steps:
\begin{align}
\sum_{i\in L}\frac{\bar T_i^{x^n}(t_2) -\bar T_i^{x^n}(t_1)}{\pi_i(1)(t_2-t_1)}\label{eq:T1}
&=
\sum_{i\in L}\frac{\left|A_i^{x^n}\cap (t_1^n,t_2^n]\right|}{\pi_i(1)(t_2^n-t_1^n)}\label{eq:T2}
\\
& \approx
\sum_{i\in L}\frac{\left|A_i^{x^n}\cap (\tau_n,t_2^n]\right|}{\pi_i(1)(t_2^n-\tau_n)}
\label{eq:T3}
\\
& =
\sum_{i\in L}\frac{\left|
\{t\in (\tau_n,t_2^n]:
C^{x^n}(t)=i,
E_i^{x^n}(t)=1,
Q_i^{x^n}(t)>0    
\}
\right|}{\pi_i(1)(t_2^n-\tau_n)}
\label{eq:T1}
\\
& \approx
\sum_{i\in L}\frac{\left|
\{t\in (\tau_n,t_2^n]:
C^{x^n}(t)=i,
E_i^{x^n}(t)=1    
\}
\right|}{\pi_i(1)(t_2^n-\tau_n)}
\\ \label{eqn:pata}
& \approx
\sum_{i\in L}\frac{\left|
\{t\in (\tau_n,t_2^n]:
C^{x^n}(t)=i
\}
\right|}{t_2^n-\tau_n}\xrightarrow[n\to\infty]{}1.
\end{align}
The symbol $\approx$ means `equal with an error that vanishes as $n\to\infty$'.
The first and second identities ((\ref{eq:T1}) and (\ref{eq:T3})) are respectively because of the definitions of $\bar T^{x^n}_i$ and $A_i^{x^n}$.
The first approximation (\ref{eq:T2}) is due to \eqref{eqn:cebolla}.
The validity of \eqref{eqn:ladrillo} implies that its microscopic version
\begin{align}\label{eqn:cinta}
\min_{i\in L}\bar Q^{x^n}_i(t)>
\max_{i\notin L}\bar Q^{x^n}_i(t)\quad \forall t\in [t_1,t_2]
\end{align}
holds for $n$ large enough,
which explains the second approximation
in which we simply neglect the condition of having a positive queue length.
The convergence to $1$ stated in the last line also uses  the asymptotic validity of \eqref{eqn:cinta}.
Indeed, for large values of $n$,
the SLC policy will dedicate service only to the queues in $L$ during the time interval $(\tau_n,t_2^n]$ because, on the one hand, the leading queue is necessarily in $L$ and,
on the other hand, the queue under service is connected at decision epochs as explained before.
We recall that the last approximation was already explained towards the end of Section \ref{sec:TFLtwo}.


While starting with $x^n$ as initial condition,
let $([a^n(m),b^n(m)))_{m\in\NN}$ be the sequence of intervals after $t_1\|q^n\|$ inside which all the environments are connected.
We now prove convergence \eqref{eqn:enchufla}.
Let $A_i^n(m)$ be the event defined as follows: in the time interval $[a^n(m),b^n(m))$, there are exactly one $\lambda_i$-mark and exactly one $\mu_i$-mark, and the $\lambda_i$-mark comes before the $\mu_i$-mark.
Let $A^n(m)=\bigcap_{i\in[K]}A_i^n(m)$.
For $n$ large enough such that \eqref{eqn:cinta} holds, the occurrence of $A^n(m)$ for an $m$ such that $b^n(m)\le t_2\|q^n\|$ implies that $C^{x^n}(b^n(m))\in L$.
Indeed, under these conditions, we can guarantee that a decision epoch occurs inside $[a^n(m),b^n(m))$, epoch at which we will pass to serve the maximum queue (that is in $L$) because all the queues are connected.
In conclusion, if $A^n(m)$ occurs then $\tau^n\le b^n(m)$, or, in other terms,
\begin{align}
\tau^n\le \min\{b^n(m):m\in\NN,b^n(m)\le t_2\|q^n\|,A^n(m)\text{ occurs}\}.
\end{align}
We can conclude from the fact that the occurrences of the events $\{A^n(m):m\in\NN\}$ represent independent trials of positive probability.

\subsection{Proof of (\ref{eqn:prisma})}
\label{app:g}

It follows from the following steps:
\begin{align}\label{eqn:tarta}
\sum_{i\in L}\frac{\bar T_i^{x^n}(t_2) -\bar T_i^{x^n}(t_1)}{\theta_i(\gamma)(t_2-t_1)}
& \approx
\sum_{i\in L}\frac{\left|
\{t\in (t_1^n,t_2^n]:
C^{x^n}(t)=i,
E_i^{x^n}(t)=1    
\}
\right|}{\theta_i(\gamma)(t_2^n-t_1^n)}
\\ 
& \approx
\sum_{i\in L}\frac{\left|
\{t\in (t_1^n,t_2^n]:
C^{x^n}(t)=i
\}
\right|}{t_2^n-t_1^n}\xrightarrow[n\to\infty]{}1-\pi_L^0.
\label{eq:rate}
\end{align}
The first approximation follows as 
in Section \ref{sct:calor}
but without the necessity of introducing the stopping time $\tau_n$.
We highlight that \eqref{eqn:ladrillo} is required in this step.
The second approximation and the convergence were already explained in Section \ref{sec:TFLthree}.

\subsection{Test verification for $\pol=\mr I$}

\label{subsection:pol_1}

We will verify the TFL for 
\begin{align}
\delta=\frac{\min\{1-\pi_L^{0}-\sum_{i\in L}\rho_i:L\subseteq[K],L \neq\emptyset\}}
{\sum_{i=1}^K\mu_i^{-1}},
\end{align}
which is a positive quantity due to condition \eqref{eqMSRI}.
Let $t_1<t_2$ satisfying \eqref{eqn:ladrillo}, and call $L=L(t_1)=L(t_2)$.
As in the other settings, we have
\begin{align}\label{eqn:mouse}
\frac{M(t_2)-M(t_2)}{t_2-t_1}\sum_{i\in L}\mu_i^{-1}
=\sum_{i\in L}\rho_i
-\sum_{i\in L} \frac{T_i(t_2)-T_i(t_1)}{t_2-t_1},
\end{align}
and the proof is now reduced to proving that
\begin{align}\label{eq:jabalina}
\sum_{i\in L} \frac{T_i(t_2)-T_i(t_1)}{t_2-t_1}
=1-\pi_L^{0}.
\end{align}
As before, condition \eqref{eqn:ladrillo}
let us obtain this identity through the following steps:
\begin{align}
\sum_{i\in L}\frac{\bar T^{x^n}_i(t_2)-\bar T^{x^n}_i(t_1)}{\mu_i(t_2-t_1)}
&\approx \frac{\sum_{i\in L}\left|\left\{t\in(t_1^n,t_2^n]:C^{x^n}(t)=i,E^{x^n}_i(t)=1\right\}\right|}{t_2^n-t_1^n}
\\ \label{eqn:truco}
&\approx\frac{\left|\left\{t\in(t_1^n,t_2^n]:E^{x^n}_i(t)=1\text{ for some }i\in L\right\}\right|}{t_2^n-t_1^n}
\xrightarrow[n\to\infty]{}1-\pi_L^{0}.
\end{align}

\section{Additional proofs}
\label{sec:addproofs}

\subsection{Proof of Corollary~\ref{eq:corredIII}}

We first note that  vertices characterizing both $\MSR{I}$ and $\MSR{III}^\PPP(\gamma)$ are related. Indeed, if
$(x_1,x_2,\ldots,x_K)$ is a vertex of $\MSR{I}$, it then follows that $(\theta_1(\gamma) x_1, \theta_2(\gamma) x_2,\ldots, \theta_K(\gamma) x_K)$ is a vertex of $\MSR{III}^\PPP(\gamma)$. 
To see this it suffices to observe that with the change of variable $\tilde \rho_i =\rho_i/\theta_i(\gamma)$, the inequalities of Setting III, see~Equation~(\ref{eqMSRIII}),  coincide with those of Setting I, see~Equation~(\ref{eqMSRI}).

This now allows us to relate their volumes. The volume is a measure of the "size" of the convex set and if we scale each dimension by a factor $\theta_i(\gamma)$,
the volume scales by the product $\theta_1(\gamma) \times \cdots \times \theta_K(\gamma)$. In other words, $\VOL{I}$~and~$\VOL{III}(\gamma)$ satisfy the relation $\VOL{III}(\gamma) = \theta_1(\gamma) \theta_2(\gamma) \cdots  \theta_K(\gamma)  \VOL{I}$, which directly yields the expression for $\SR{III}(\gamma)$ in Equation~\eqref{eq:redIII}.

\subsection{Proof of Corollary \ref{optigamma}}

We can follow the exact same steps of the proof of stability of the system without delays. In equation (\ref{eq:rate}), the limit is replaced by ${1/\gamma \over {1\over \gamma} +c} (1-\pi^{0}_L)$.


\subsection{Proof of Proposition \ref{prop:strict}}

The proof is very similar to the test verification for $\pol=\rm{III}$ given in Section \ref{subsection:pol_2_and_3}.
With the new policy, 
the RHS of \eqref{eqn:tarta}
need to be decomposed as
\begin{align}
&\sum_{i\in L}\frac{\left|\left\{t\in(t_1^n,t_2^n]:C^{x^n}(t)=i,E^{x^n}_i(t)=1,t\in B_n\right\}\right|}{\theta_i(\gamma)(t_2^n-t_1^n)}
\\
&\quad +
\sum_{i\in L}\frac{\left|\left\{t\in(t_1^n,t_2^n]:C^{x^n}(t)=i,E^{x^n}_i(t)=1,t\notin B_n\right\}\right|}{\theta_i(\gamma)(t_2^n-t_1^n)}.
\end{align}
Here $B_n$ is the union of the intervals delimited by consecutive $\gamma$-marks for which all the queues are disconnected at the beginning of the interval.
The second sum is controlled as we did before, with the only difference that the concerning tailor-made Markov process only considers the time outside $B_n$.
To control the first sum,
we need to define a similar process but that only records the time inside $B_n$.
For the same reasons than before, the asymptotic proportion of time that such a new process is connected is
\begin{align}
\phi_i(\gamma)=\frac{\lambda_i'}{\lambda_i'+\mu_i'+\gamma}.
\end{align}
If we define
\begin{align}
\varepsilon(\gamma)=\min_{i\in [K]}\frac{\phi_i(\gamma)}{\theta_i(\gamma)}
=
\min_{i\in [K]}\frac{\lambda_i'}{\gamma+\lambda_i},
\end{align}
then the second sum is asymptotically larger or equal than $\varepsilon(\gamma)\pi_L^0$.
Repeating the same steps than before, we obtain that the region
\begin{align}
\sum_{i\in L}\frac{\rho_i}{\theta_i(\gamma)}<1-\pi_L^0(1-\varepsilon(\gamma)) \quad \forall L\subseteq [K],L\neq\emptyset,
\end{align}
which of course strictly contains $\MSR{III}^\PPP(\gamma)$,
is contained in the stability region of the policy at issue.

To prove the second inequality,
we observe that
\begin{align}
\varepsilon'(\gamma)=\min_{i\in [K]}\phi_i(\gamma)
\end{align}
gives a lower bound for the proportion of capacity wasted when the $\gamma$-decision epochs are imposed, whatever the policy is and whichever queue is receiving service.
Fix a queue $i$, take $\rho_i$ such that
\begin{align}
1-\pi_i(0)-\varepsilon'(\gamma)<\rho_i<1-\pi_i(0),
\end{align}
and  chose the rest of the $\rho_j$'s to be inside $\MSR{I}$.
Under these choices,
we are of course inside the maximal stability region for Setting I,
but at the same time tasks in queue $i$ accumulate in Setting III,
and hence we are outside $\MSR{III}(\gamma)$ as desired.

\end{document}